\documentclass[11pt,a4paper]{amsart}

\newcommand{\figdir}{}

\usepackage[margin=3cm]{geometry}
\usepackage{amsmath,amsfonts,amssymb,amsthm}
\usepackage{graphicx}
\usepackage{multirow}
\usepackage{subfigure}
\usepackage{float}
\usepackage{rotating}

 \usepackage{xspace}

 \usepackage{epsfig}

 \def\R{\mathbb R}
\newcommand{\pdpd}[2]{\frac{\partial #1}{\partial #2}}

\begin{document}
\title{A multigrid method for the Helmholtz equation with optimized coarse grid corrections}
\author{Christiaan C. Stolk}
\address{Korteweg-de Vries Institute for Mathematics, University of Amsterdam,
Science Park 904, 1098 XH Amsterdam}
\author{Mostak Ahmed}
\address{Department of Mathematics, Jagannath University, Dhaka-1100, Bangladesh}
\author{Samir Kumar Bhowmik}
\address{Department of Mathematics, University of Dhaka, Dhaka-1000, Bangladesh and
Department of Mathematics and Statistics, College of Science, Al Imam Mohammad Ibn Saud Islamic University, P.O. BOX 90950, 11623 Riyadh, KSA}

\maketitle

\begin{abstract}
We study the convergence of multigrid schemes for the Helmholtz equation, focusing in particular on the choice of the {\em coarse scale operators}. Let $G_{\rm c}$ denote the number of points per wavelength at the coarse level. If the coarse scale solutions are to approximate the true solutions, then the oscillatory nature of the solutions implies the requirement $G_{\rm c} > 2$. However, in examples the requirement is more like $G_{\rm c} \gtrsim 10$, in a trade-off involving also the amount of damping present and the number of multigrid iterations. We conjecture that this is caused by the difference in phase speeds between the coarse and fine scale operators. Standard 5-point finite differences in 2-D are our first example. A new coarse scale 9-point operator is constructed to match the fine scale phase speeds. We then compare phase speeds and multigrid performance of standard schemes with a scheme using the new operator. The required $G_{\rm c}$ is reduced from about 10 to about 3.5, with less damping present so that waves propagate over $>$ 100 wavelengths in the new scheme. Next we consider extensions of the method to more general cases. In 3-D comparable results are obtained with standard 7-point differences and optimized 27-point coarse grid operators, leading to an order of magnitude reduction in the number of unknowns for the coarsest scale linear system. Finally we show how to include PML boundary layers, using a regular grid finite element method. Matching coarse scale operators can easily be constructed for other discretizations. The method is therefore potentially useful for a large class of discretized high-frequency Helmholtz equations.
\end{abstract}



\markboth{C.C. STOLK, M. AHMED AND S.K. BHOWMIK}{MULTIGRID FOR THE HELMHOLTZ OPERATOR}

\section{Introduction}

Large scale Helmholtz problems are notoriously difficult to solve. In particular in the high-frequency limit, when both the grid size and the frequency become large, extensive research is going on, as classical methods perform poorly and more recent methods remain costly
in the sense that they require many iterations to converge or are memory-intensive.  High-frequency Helmholtz problems have applications in various simulation problems and in inverse problems, e.g.\ in exploration seismology and acoustic scattering. 

Multigrid methods for the Helmholtz equation have been investigated by a number of authors, see e.g.\
\cite{ElmanErnstOLeary2001,ErlanggaOosterleeVuik2006,HaberMacLachlan2011,CalandraEtAl2013,
ErnstGanderPreprint2013}. Restricting to finite difference problems, a key parameter is the minimum number of grid points per wavelength in the coarsest grid, denoted here by $G_{\rm c}$. If the coarse problem solutions are to approximate the solutions of the original problem, the oscillatory nature of the solutions leads to the requirement $G_{\rm c} > 2$. In this parameter regime fast convergence is possible. However, the coarse grid problem can remain large, especially since for existing methods the requirement is more like $G_{\rm c} \gtrsim 10$. This regime is therefore mostly of interest when the discretization has many points ($\gtrsim 20$) per wavelength, due e.g.\ to subwavelength detail in the coefficient or the right-hand side. 
In the high-frequency case smaller values of $G_{\rm c}$ are needed to make multigrid useful. 

In a second parameter regime more multigrid levels are added and the parameter $G_{\rm c}$ is chosen $< 1$. In this regime the cost per iteration can be small, but the number of iterations is usually large. A popular method of this kind is the shifted Laplacian methods \cite{ErlanggaOosterleeVuik2006}, in which a multigrid cycle for a modified Helmholtz equation is used as a preconditioner for an iterative method like BiCGSTAB. See also \cite{ElmanErnstOLeary2001}.

Among other methods an interesting class is formed by the double sweeping domain decomposition methods, like the moving PML method of \cite{EngquistYing2011} or the double sweeping method of \cite{Stolk2013}. These methods appear to have the best, near-linear scaling for the cost per solve, but have the disadvantage of a large memory use, see also \cite{PoulsonEtAl2012_Preprint}. Other methods use e.g.\ incomplete factorizations and matrix compression techniques or a combination of techniques \cite{BollhoferGroteSchenk2009,ChenChengWu2012,WangDeHoopXia2011,CalandraEtAl2013}. 
In all cases it is important to distinguish the behavior in non-resonant cases vs.\ resonant cases. In variable coefficient media with resonances the performance of the iterative methods discussed above tends to deteriorate strongly.

In this paper we study multigrid methods in the first regime, aiming at fast convergence, say $\lesssim 20$ iterations for reduction of the residual by $10^{-6}$. Our contribution centers on the choice of the {\em coarse scale operators}. We describe a simple criterion for the choice of the coarse scale operator. Using new, optimized coarse scale operators and carefully selected smoother parameters, we will obtain a reduction of $G_{\rm c}$ from about 10 to about 3.5 and still have fast convergence. As examples we study the standard 5-point and 7-point stencils (in 2-D and 3-D resp.) for the fine scale operator, and a regular grid finite element method for which we show how to include PML boundary layers. 

The main consequence is that the two-grid method can now be used in high-frequency Helmholtz problems. It becomes a general method to reduce the number of degrees of freedom in a Helmholtz problem, see also the discussion in Section~\ref{sec:discussion}. In situations where the multigrid method was already useful, e.g.\ when subwavelength detail is present, the number of degrees of freedom in the expensive coarse scale problem can potentially be reduced by a factor of $(10 / 3.5)^3 \approx 23$ in 3-D.

\medskip

We now introduce the setup in more detail. The Helmholtz equation reads
\begin{equation} \label{helm_with_i}
  L u 
  \stackrel{\rm def}{=}
  - \Delta u - ((1+\alpha i)k(x) )^2 u = f ,
\mbox{\quad in }\Omega\subset\R^n
\end{equation}
where $\Delta$ is the Laplacian and 
$\Omega$ is a rectangular block. We assume Dirichlet boundary conditions at the boundary $\partial \Omega$. PML layers will be present in some of the examples, they will be discussed
in section~\ref{sec:further_numerics}.
The parameter $k$ in general depends on $x \in \Omega$.
Here $\alpha$ is a parameter for the damping, that will mostly be independent of $x$%
\footnote{Other authors sometimes use $-\Delta - (1+i\alpha)k^2$, i.e.\
with the factor $(1+i\alpha)$ outside the square. The sign of $i\alpha$
is related to our choice of temporal Fourier transform  
$f(t) = e^{-i \omega t} \hat{f}(\omega)$.}.
The corresponding undamped operator will be denoted by
\begin{equation}
  H = - \Delta - k^2 .
\end{equation}
The imaginary contribution $i\alpha k$ leads to exponentially decaying solutions. For example, in 1-D for constant $k$ there are solutions
\begin{equation}
  e^{i k x - \alpha k x} .
\end{equation}
This limits the propagation distance of the waves. This distance, measured by the number of wavelengths for the amplitude to be reduced by a factor 10, is given by
\begin{equation}
  D(\alpha)= \frac{\log(10)}{2 \pi \alpha } .
\end{equation}
Values of $\alpha$ will range from $1.25 \cdot 10^{-3}$ to $0.02$, small enough for applications.
In presence of PML layers $\alpha$ will be set to zero. The experiments will show that, for the optimized coarse grid methods good convergence can already start around $\alpha=1.25 \cdot 10^{-3}$  or at $D(\alpha) \approx 290$ wave lengths assuming around 3.5 points per wavelength in the coarsest grid. 

Multigrid methods consist of several components \cite{TrottenbergOosterleeSchueller2001}. For example, a two-grid method consists of the following steps: presmoothing using a relaxation method; restriction of the residual to the coarse grid; solving a coarse grid equation; prolongation of the solution to the fine grid; suppressing remaining errors using postsmoothing, using again a relaxation method. Multigrid methods are used by themselves or as a preconditioner for a Krylov subspace solver. The latter option will be adopted in this paper. 

Our study focusses on two parts of the multigrid method: The coarse grid operator and the smoother. Concerning the coarse grid operator our claim is that it should have the same {\em phase speed} or numerical dispersion as the fine scale problem. Indeed over large distances the phase errors lead to large differences between the approximate, coarse grid solution and the true solution. Alternatively, from Fourier analysis of multigrid methods one can argue that the {\em inverse of the symbols} for the coarse scale and original operators should match, which in turn implies that the zero-sets and hence the phase speeds should match. Standard choices for the coarse scale operators are the Galerkin approximation, or to use the same discretization scheme as the fine scale operator. They lead to sizeable phase speed differences. As shown in section~\ref{sec:phase_speeds}, these differences can be sharply reduced by using finite difference schemes with optimized coefficients instead.

When designing the smoother there are many choices to be made, concerning the method and the parameter values. In this work we consider the standard SOR and $\omega$-Jacobi methods. The parameters involved are then the relaxation parameter, denoted by $\omega_{\rm S}$, and the number of pre- and postsmoothing steps $(\nu_1,\nu_2)$. The effect of the different choices is studied mainly using local Fourier analysis (LFA) in section~\ref{sec:MG_Fourier_analysis}. The results show that with the optimized coarse grid operators and standard smoothers it is indeed possible to have a converging method for small $G_{\rm c}$ and $\alpha$, provided that suitable parameters are chosen. For small $G_{\rm c}$ and $\alpha$ the convergence depends quite sensitively on these parameters. 


The first example we take is the two-grid method with standard 5-point finite differences in 2-D. 
We first compute the optimized coefficients, and compare the phase speeds associated with optimized and standard coarse grid operators. Then the two-grid method is studied for different choices of the parameters in two ways, first the convergence factors are computed, and finally the convergence in numerical experiments is studied. The multigrid convergence is good or poor, precisely when the phase speeds match well or poorly, respectively. The two-grid convergence factors show this for constant $k$. For spatially varying $k$ this is verified using numerical experiments.

Next we consider more general settings than 5-point finite differences. For the 3-D problem we construct optimized coarse scale operators for the standard 7-point finite difference scheme. Numerical tests show similar results as for the 2-D case, with good convergence for $G_{\rm c} \gtrsim 3.5$. In practice problems on rectangular domains often occur in combination with absorbing boundary conditions or absorbing boundary layers, such as PML layers \cite{Berenger1994,ChewWeedon1994}. These then provide damping in the equation, so that we set $\alpha = 0$. In numerical experiments, straightforward insertion of a PML layer was observed to lead to very poor convergence or no convergence at all in the situations we are interested in ($G_{\rm c} \gtrsim 3.5$). To address this we consider an adapted coarsening strategy. In the PML layer no coarsening takes place in the direction normal to the boundary. This is the direction in which a fast variation of the coefficient $\sigma$ used in the PML takes place. Because this leads to (partly) irregular grids, it is natural to consider this in the context of finite elements. Optimized coarse grid operators for a first order rectangular finite element discretization are computed, and we show in 2-D examples that this again results in good convergence of the multigrid method for $G_{\rm c} \gtrsim 3.5$, and propapagation distances of up to 200 wavelengths. 

The setup of the paper is as follows. The next section focusses on the phase speeds. Here we construct the new, optimized finite difference operators for the various cases and compare the phase speed errors in the standard and the new coarse grid operators. Then in section~\ref{sec:multigrid_methods} we describe the multigrid methods used in this study. In section~\ref{sec:MG_Fourier_analysis} we present Fourier analysis of the two-grid methods. Section~\ref{sec:numerical_results} describes the results of numerical simulations: First the standard and the new methods are compared for 2-D finite differences, then the extension to multigrid, to 3-D and to problems with PML boundary layers is discussed. We end the paper with a brief discussion section.



\section{Phase speeds and optimized coarse scale operators}
\label{sec:phase_speeds}

We first recall the notions of dispersion relation and phase speed for constant coefficient linear time dependent partial differential equations \cite{Trefethen1982}. If $P$ is such an operator, and 
\begin{equation}
  \sigma_P(\xi,\omega)  = e^{ - i \xi \cdot x + i \omega t}
  ( P e^{ i \xi \cdot x - i \omega t} ) , 
\end{equation}
denotes its symbol, then the dispersion relation is the set of $(\xi,\omega)$ where $\sigma_P(\xi,\omega) = 0$.
For the Helmholtz equation, the symbol is a function of the spatial wave number $\xi$, with parameter $\omega$
\begin{equation} \label{eq:def_Helm_symbol}
  \sigma_H(\xi;\omega) = e^{ - i \xi \cdot x} ( H e^{ i \xi \cdot x} ) 
\end{equation}
and the dispersion relation is understood to be the $\omega$-dependent set
\begin{equation} \label{eq:def_Helm_disprel}
  \{ \xi \in S \, | \, \sigma_H(\xi;\omega) = 0  \} ,
\end{equation}
where for the continuous operator $S = \R^n$, and it is assumed that 
$\alpha = 0$. For $\xi$ such that $\sigma_H(\xi;\omega) = 0$, the number 
$\frac{\omega}{\| \xi \|}$ is called the phase speed $v_{\rm ph}$ associated with a plane wave solution. For the continuous operator $v_{\rm ph}$ is constant equal to $c$.

For finite difference discretizations of the Helmholtz operator
the definitions (\ref{eq:def_Helm_symbol}) and 
(\ref{eq:def_Helm_disprel}) remain valid except that
$S$ is given by the fundamental domain
$S = [-\pi/h,\pi/h]^n$.
In the typical case that $\sigma_H$ is increasing along half lines from the 
origin, the phase speed is a function of angle. We will compute the dimensionless phase speed $\frac{v_{\rm ph}}{c}$ as a function of another dimensionless quantity, the number of points per wavelength $G$, or its inverse $1/G$ ($G$ is related to the dimensionless quantity $k h$ by $k h = \frac{2 \pi}{G}$). When multiple grid levels are involved $G_{\rm c}$ refers to the coarse level.

This discussion will mostly involve only two scales, a fine scale and a coarse scale with double the grid parameter. In multigrid with more than two levels, the additional levels are assumed to be increasingly fine, because the coarsest level is restricted by the wave length. In general we expect that the phase speed difference between the two coarsest levels is the most important, since for finer grids these differences automatically become smaller as the discretization becomes a better approximation of the true operator.

Next we first show the behavior of the phase speeds for some well known finite difference schemes and then construct a finite difference method that is optimized so that its dispersion relation matches that of the standard 5-point method. After this we consider the 7-point operator in 3-D and the finite element method used with PML layers in 2-D.

\subsection{Phase speeds of some well known finite difference schemes}

\begin{figure}[thf]
\begin{minipage}[t]{6cm}
\centering
(a)\\
\includegraphics[height=44mm]{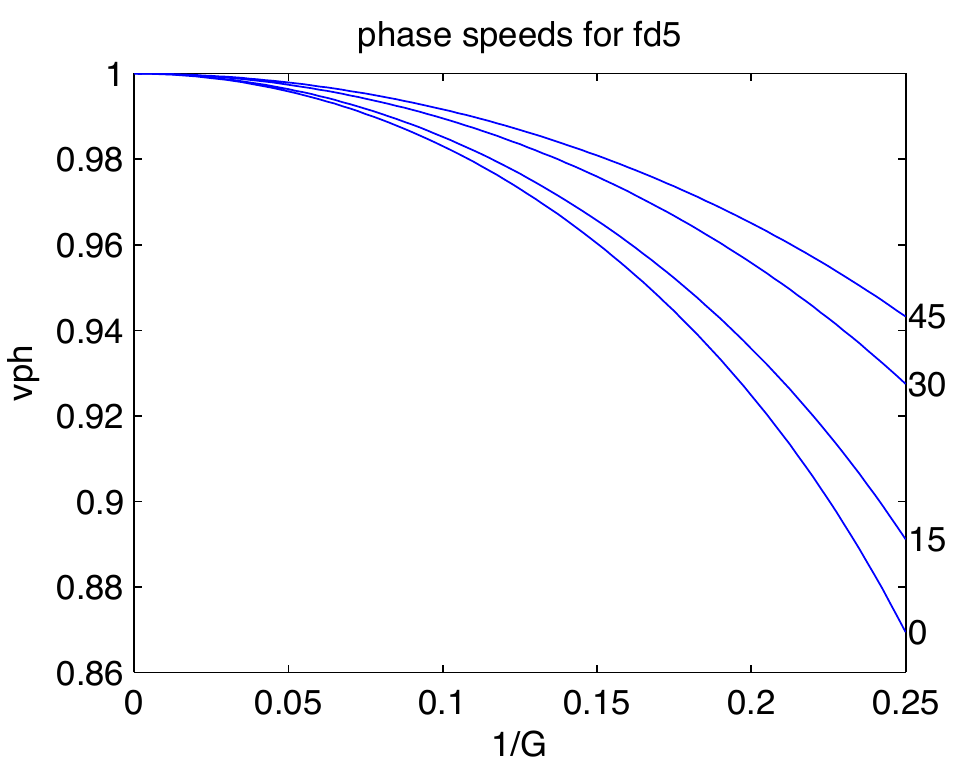}
\end{minipage}
\hspace*{2mm}
\begin{minipage}[t]{6cm}
\centering
(b)\\
\includegraphics[height=44mm]{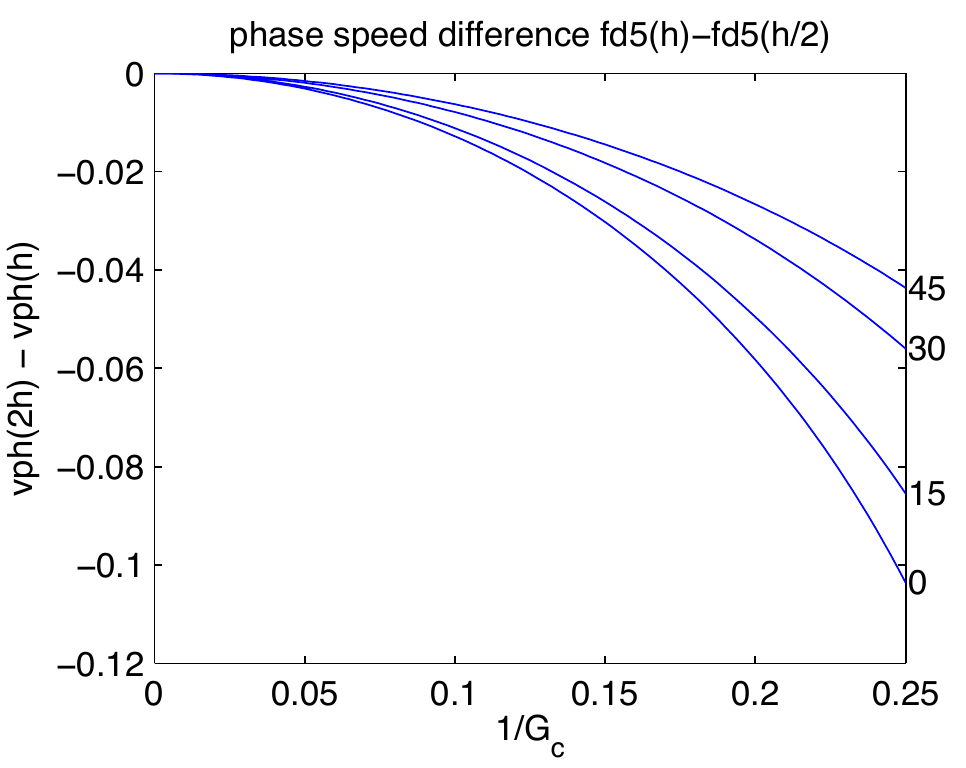}
\end{minipage}\\
\begin{minipage}[t]{6cm}
\centering
(c)\\
\includegraphics[height=44mm]{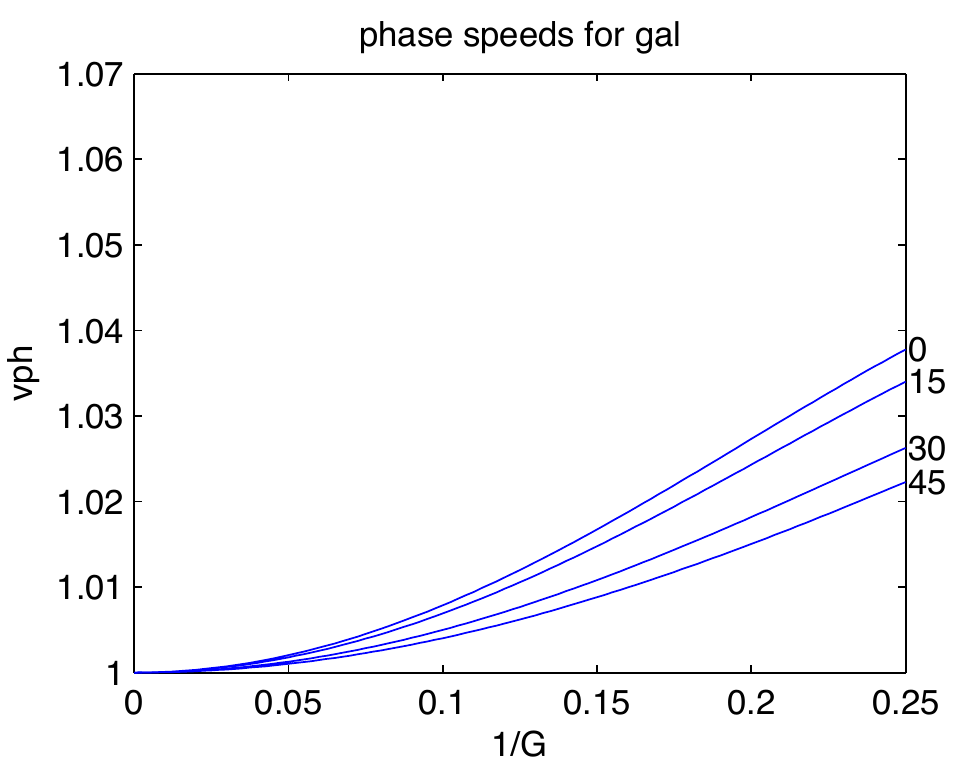}
\end{minipage}
\hspace*{2mm}
\begin{minipage}[t]{6cm}
\centering
(d)\\
\includegraphics[height=44mm]{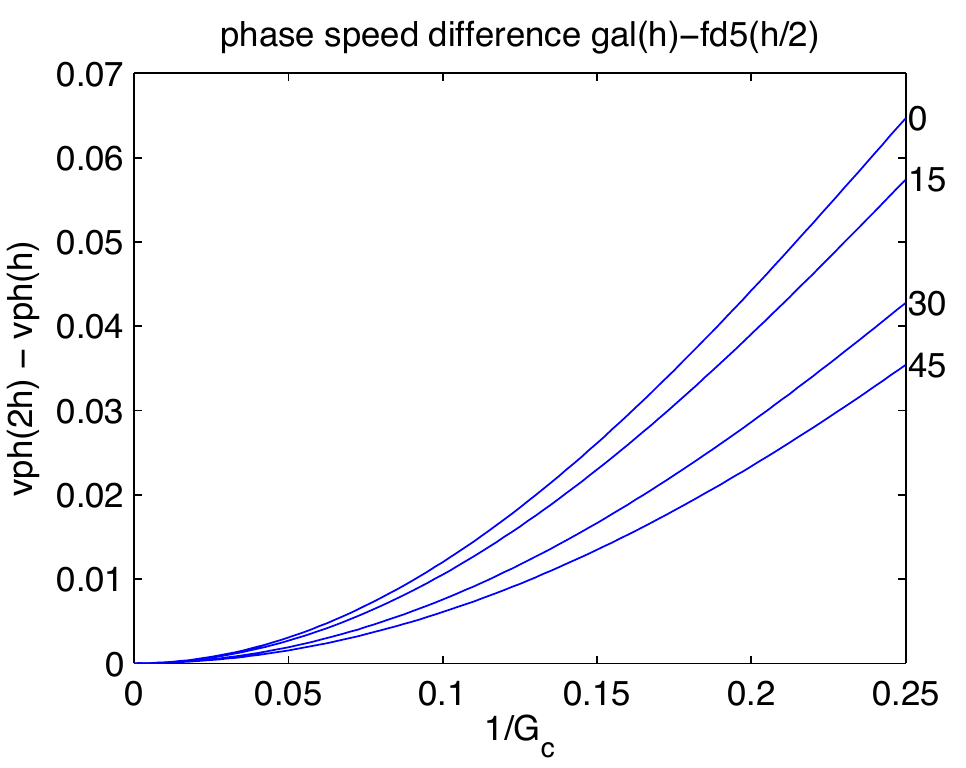}
\end{minipage}\\
\begin{minipage}[t]{6cm}
\centering
(e)\\
\includegraphics[height=44mm]{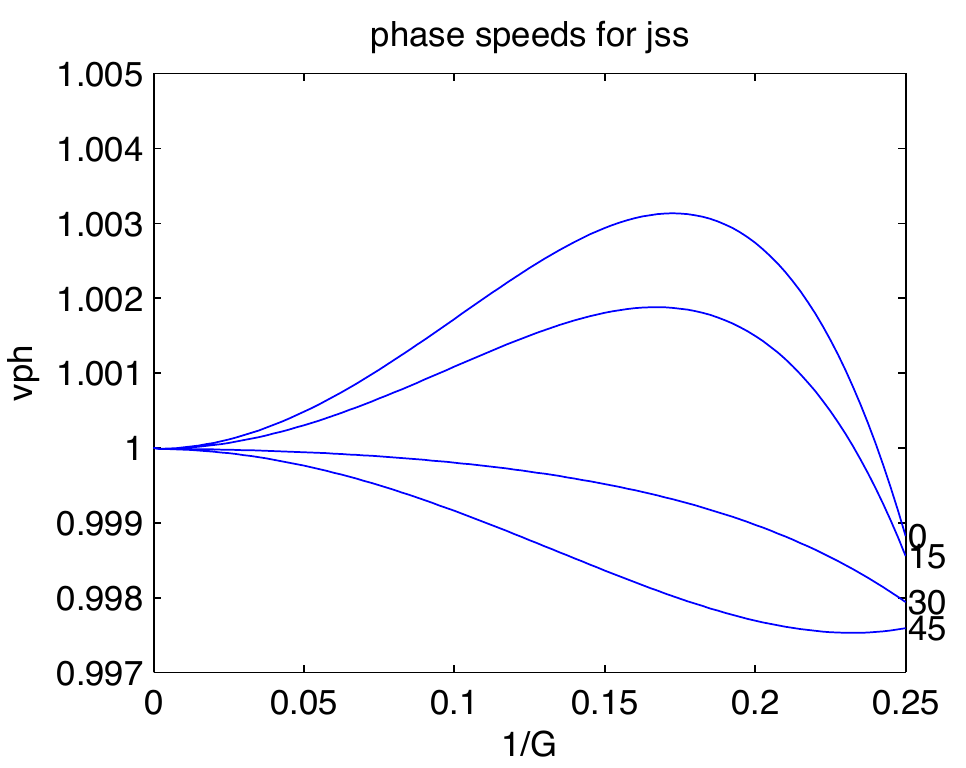}
\end{minipage}
\hspace*{2mm}
\begin{minipage}[t]{6cm}
\centering
(f)\\
\includegraphics[height=44mm]{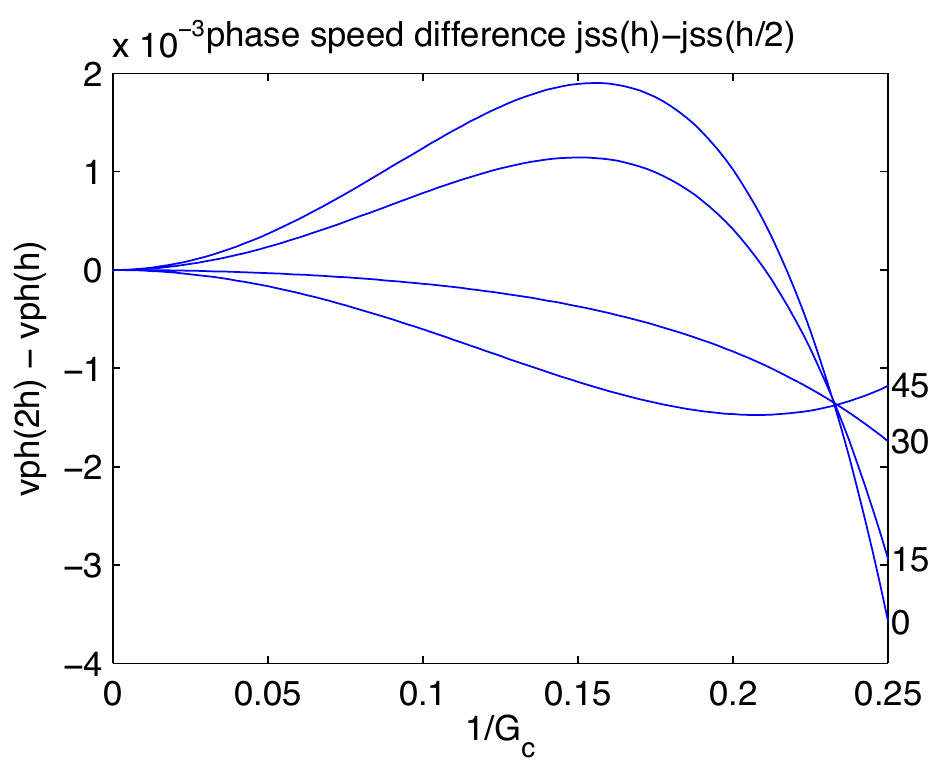}
\end{minipage}

\caption{Dimensionless phase speed curves and difference between fine and coarse scale phase speeds; (a), (c) and (e): phase speed for the fd5, Galerkin and JSS method; (b), (d) and (f): phase speed differences $v^{\rm fd5}_{h} - v^{\rm fd5}_{h/2}$, $v^{\rm gal}_{h} - v^{\rm fd5}_{h/2}$, and 
$v^{\rm JSS}_{h} - v^{\rm JSS}_{h/2}$, with angles $0^\circ$, $15^\circ$, $30^\circ$, $45^\circ$}
\label{fig:disprel1}
\end{figure}

The standard 5-point finite difference discretization (fd5) of the 2-D Helmholtz operator is given by
\begin{equation}
  (H^{\rm fd5}_h u)_{i,j} = h^{-2} ( 4 u_{i,j} - u_{i-1,j} - u_{i+1,j} - u_{i,j-1} - u_{i,j+1} ) - k^2 u_{i,j} .
\end{equation}
Its symbol is easily shown to be 
\begin{equation}
  \sigma^{\rm fd5}_h (\xi) = h^{-2} ( 4 - 2 \cos(h \xi_1) - 2 \cos(h \xi_2) ) - k^2 ,
\end{equation}
where $\xi = (\xi_1,\xi_2)$ denotes the wave vector.  The phase speed as a function of angle is easily computed using a root finding algorithm and shown in Figure~\ref{fig:disprel1}(a). The phase speeds as a function of $|\xi|$ are given for 4 angles $0^\circ$, $15^\circ$, $30^\circ$ and
$45^\circ$.

When such a scheme is used in a multigrid method, say a two-grid method for the purpose of this argument, standard choices for the coarse scale operator are to use the same discretization, or to use a Galerkin discretization. The comparison between the phase speeds of a coarse scale operator $H^{\rm fd5}_h$ and a fine scale operator $H^{\rm fd5}_{h/2}$ is given in Figure~\ref{fig:disprel1}(b). The Galerkin method, using ``full weighting'' restriction and prolongation operators \cite{TrottenbergOosterleeSchueller2001} is easily shown to have symbol
\begin{equation}
\begin{split}
\sigma^{\rm gal}_h(\xi) = 
{}& 
3h^{-2} - \frac{9}{16} k^2  
+ ( - h^{-2} - \frac{3}{16} k^2  ) (\cos(\xi_1 h) + \cos(\xi_2 h) )
\\
{}& 
+ ( - \frac{1}{2} h^{-2} - \frac{1}{32} k^2)
 ( \cos( (\xi_1+\xi_2)h) + \cos( (\xi_1-\xi_2)h) ) .
\end{split}
\end{equation}
The phase speeds and the phase speed difference between a coarse scale Galerkin and a fine scale fd5 method are shown in Figure~\ref{fig:disprel1}(c) and (d).  The phase speed differences are on the order of 0.01 or 0.02 for $G_{\rm c} = 8$ and larger for $G_{\rm c}$ smaller.

We contrast this with the optimized finite difference method of Jo, Shin and Suh \cite{JoShinSuh1996}, sometimes called the mixed grid operator. This operator, here called JSS operator for brevity, is given by
\begin{equation} \label{eq:JSS_discrete_op}
\begin{split}
  (H^{\rm jss}_h u)_{i,j} 
  = {}& 
( (2 + 2a) h^{-2} - c k^2 ) u_{i,j}
+ ( - a h^{-2} - d k^2 ) ( u_{i+1,j} + u_{i-1,j} + u_{i,j+1} + u_{i,j-1} )
\\
{}& + ( - \frac{1-a}{2}h^{-2} - \frac{(1-c-4d)}{4} k^2 )
 ( u_{i+1,j+1} + u_{i-1,j+1} + u_{i+1,j-1} + u_{i-1,j-1} )
\end{split}
\end{equation}
where $a = 0.5461$, $c = 0.6248$ and $d = 0.9381 \cdot 10^{-1}$.
The symbol of this operator hence equals
\begin{equation}
\begin{split}
  \sigma^{\rm jss}_h(\xi) = 
{}&
 ( h^{-2} (2 + 2a) - k^2 c ) 
+ 2 ( - h^{-2} a - k^2 d ) (\cos(\xi_1 h) + \cos(\xi_2 h) )
\\
{}& + 2 ( -h^{-2} \frac{1-a}{2} - \frac{(1-c-4d)}{4} k^2 )
 ( \cos( (\xi_1+\xi_2)h) + \cos( (\xi_1-\xi_2)h) ) .
\end{split}
\end{equation}
The phase speed and relative phase speed difference between
$H^{\rm jss}_{2h}$ and $H^{\rm jss}_h$ are given in Figure~\ref{fig:disprel1}(e) and (f).
The phase speed difference between coarse and fine level operators, both using the JSS method, is reduced to around $2.5 \cdot 10^{-3}$ for $G_{\rm c}$ down to 4, i.e.\ a very substantial improvement both in $G_{\rm c}$ and in the size of the phase speed difference.

The operators $L^{\rm fd5}_h$, $L^{\rm gal}_h$ and $L^{\rm jss}_h$, with $\alpha \neq 0$ are obtained simply by replacing $k^2$ by $( (1 + i \alpha ) k)^2$.

\subsection{An optimized finite difference scheme}

To obtain a finite difference method with phase speed matching that of the standard 5-point method, we will now consider certain finite difference discretizations with optimized coefficients depending on $h$, $k$ and the fine grid parameters $h_{\rm f}$. The new operators match the phase speed of the standard 5-point operator with very good accuracy, even down to $G_{\rm c}=2.5$.

For each of the two terms in $-\Delta -k^2$ we describe the discretization. Like in \cite{JoShinSuh1996} the discretization of the mass term $-k^2$ involves a symmetric 9-point stencil. This stencil depends on three coefficients $b_1,b2$ and $b_3$ as follows
\begin{equation} \label{eq:2D_mass_sten_b}
  \begin{bmatrix} b_3/4 & b_2/4 & b_3/4 \\ b_2/4 & b_1 & b_2/4 \\ b_3/4 & b_2/4 & b_3/4
  \end{bmatrix} .
\end{equation}
The coefficients satisfy $b_1+b_2+b_3=1$ (this explain the normalization in (\ref{eq:2D_mass_sten_b})) and are otherwise to be determined. The Laplacian is written as $-\Delta = -\pdpd{^2}{x_1^2} - \pdpd{^2}{x_2^2}$. The discrete second derivative $-\pdpd{^2}{x_1^2}$ is given by the tensor product of a 1-D mass matrix with stencil $\begin{bmatrix} a_2/2 \\ a_1 \\ a_2/2 \end{bmatrix}$ and the standard second order derivative, with stencil $\big[ -h^{-2} \; 2h^{-2} \; -h^{-2} \big]$. The coefficients are assumed to be such that $a_1+a_2 = 1$, otherwise they are again to be determined. The second derivative $- \pdpd{^2}{x_2^2}$ is discretized using the 90 degree rotated stencil. The discrete Helmholtz operator then reads
\begin{equation} \label{opt_fd_helm}
\begin{split}
    H^{\rm opt}_h u (x_{i,j}) 
  = {}& (4 a_1 h^{-2} - k^2 b_1) u_{i,j} 
        + ( (-a_1+a_2) h^{-2} - k^2 b_2/4) (u_{i-1,j}+u_{i+1,j}+u_{i,j-1}+u_{i,j+1})\\
    {}& + ( -a_2 h^{-2} - k^2 b_3/4 ) (u_{i-1,j-1}+u_{i-1,j+1}+u_{i+1,j-1}+u_{i+1,j+1}) .
\end{split}
\end{equation}
When the coefficients are constants, this is a new description of the class of operators considered in \cite{JoShinSuh1996}. However, to obtain a better approximation of the phasespeed the coefficients $a_1,a_2,b_1,b_2,b_3$ will be allowed to {\em depend on $k$, $h$ and the fine-grid parameter $h_{\rm f}$}.

From dimensional considerations, the coefficients $a_j$, $b_j$ in fact only depend on two parameters, namely $\frac{h_{\rm f}}{h}$ and $h k$ or $1/G_{\rm c} = \frac{h k}{2 \pi}$. In our application, $\frac{h_{\rm f}}{h}$ will be a power of $(1/2)$, depending on the number of multigrid levels considered. The optimized coefficients will be computed separately for each value of $\frac{h_{\rm f}}{h}$. In this section, denote $p= 1/G_{\rm c}$. We hence have to compute the coefficients $a_1(p), b_1(p)$ and $b_2(p)$ for $p$ in an interval $[ 0, P]$ (we have applied this with $P$ up to $0.4$, i.e.\ $G_{\rm c}$ down to 2.5 points per wavelength). We describe these functions by interpolation from a small set of data points. For this purpose we consider $n_{\rm C}$ equidistant points on $[0,P]$, given by $p_k = \frac{k-1}{n_{\rm C}-1} P$, $k=1, \ldots, n_{\rm C}$. The values $a_1(p)$ for $p\in[0,P]$ will be obtained by linear interpolation from the $n_{\rm C}$ values $a(p_k)$. The functions $b_1(p)$ and $b_2(p)$ are parameterized similarly. From the relations $a_1+a_2 = 1$ and $b_1+b_2+b_3 = 1$ the remaining two coefficients are determined. We refer to the $a_1(p_k)$ and the $b_j(p_k)$, $j=1,2$ as the control values.

We next describe the computation of phase speeds. The symbol associated to (\ref{opt_fd_helm}) is
\begin{equation}
\begin{split} 
\sigma^{\rm opt}_h(\xi_1,\xi_2) = 
{}& (4 a_1 h^{-2} - k^2 b_1) 
        + ( (-a_1+a_2) h^{-2} - k^2 b_2/4) 2 (\cos(h\xi_1)+\cos(h\xi_2))\\
    {}& + ( -a_2 h^{-2} - k^2 b_3/4 ) 2 (\cos(h(\xi_1 + \xi_2)) + \cos(h(\xi_1 - \xi_2))) .
\end{split}
\end{equation}
In order to compute the phase speed for a direction given by a unit vector $(\cos(\theta), \sin(\theta))$, the equation
\begin{equation} \label{eq:solve_phasespeed}
  \sigma^{\rm opt}_h(\xi \cos(\theta), \xi \sin(\theta) ) = 0
\end{equation}
is solved for $\xi$, denote the result by $\xi^{\rm opt}_h(\theta)$. The quotient $\frac{\xi^{\rm opt}_h(\theta)}{\omega}$ may be called the phase slowness (one over the phase speed) of the discrete operator for the direction given by $\theta \in S^1$. In numerical computations, the phase slownesses are obtained from (\ref{eq:solve_phasespeed}) using a numerical equation solver (Matlab's fsolve using the default trust-region-dogleg algorithm for when the gradient is present).
 
The objective function for the optimization of the coefficients is chosen as follows. A set of angles $\theta$ and a set of $p$ values are chosen. For each pair $(\theta,p)$ the phase slownesses for the optimized symbol and the fine scale symbol, denoted by by $\frac{\xi^{\rm opt}_h}{\omega}$ and by $\frac{\xi^{\rm fd5}_{h_{\rm f}}}{\omega}$, are computed. The relative error is defined by
\begin{equation} \label{eq:relative_phase_speed_error}
  \frac{| \xi^{\rm opt}_h(\theta) - \xi^{\rm fd5}_{h_{\rm f}}(\theta) |}{\xi^{\rm fd5}_{h_{\rm f}}(\theta)} .
\end{equation}
The objective function for the optimization is the sum of the squares of the relative errors. For the results given below we used 18 equidistant values for $\theta$ to discretize the interval $[0,\pi/2)$), and 8 times $n_{\rm C}$ equidistant values for $p$. The results depended little on the precise choice of values, once they were sufficiently large.

\begin{figure}[thf]
\centering
\includegraphics[height=44mm]{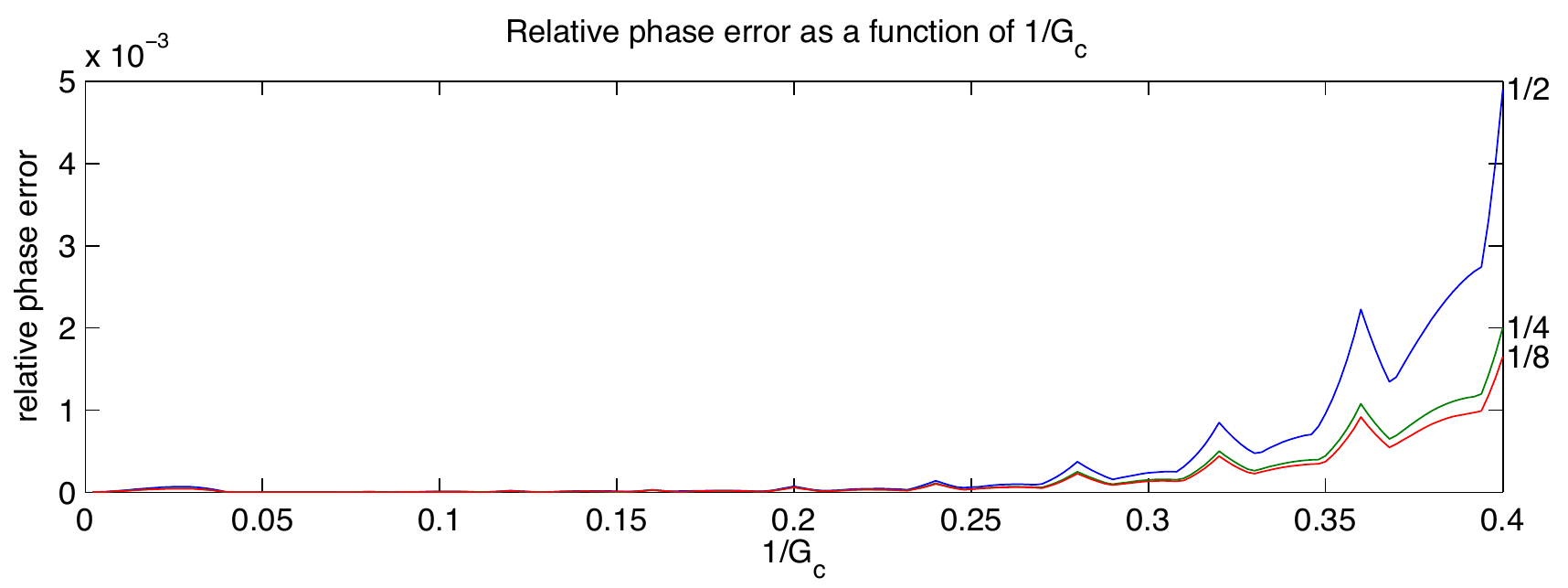}
\caption{Relative phase speed errors as defined in (\ref{eq:relative_phase_speed_error}), maximum over $\theta$ as a function of $p = 1/G_{\rm c}$, for $h_{\rm f}/h = 1/2, 1/4$ and $1/8$.}
\label{fig:err_optfd5}
\end{figure}

The objective functional is minimized as a function of the control values $a_1(p_k), b_1(p_k)$ and $b_2(p_k)$. For this purpose a constrained minimization algorithm was used, Matlab's fmincon using the interior-point algorithm. 
The behavior of the optimization algorithm depends on the parameters $P$ and $n_{\rm C}$, and on the starting values chosen. For small $P$ (less than about 0.25) the algorithm is little sensitive to the starting values, one can take for example the coefficient values from (\ref{eq:JSS_discrete_op}). By gradually increasing the value of $P$ good starting values can be found for an interval up to $P = 0.4$. The interior-point optimization algorithm then performs its job very nicely. We will use in the sequel the values found for $P=0.4$, $n_{\rm C} = 11$. They are given in Table~\ref{tab:coeff_optfd5}. The resulting errors in the relative phase speed were computed as a function of $p$, taking the maximum value as a function of $\theta$. The result is given in Figure~\ref{fig:err_optfd5}.
As can be seen from these figures, the phase speed differences are reduced to less than $2 \cdot 10^{-4}$ for $G_{\rm c}$ down to $4$ and less than $1 \cdot 10^{-3}$ for $G_{\rm c}$ down to to 3.

\begin{table}
\begin{tabular}{|r|ccc|ccc|ccc|} \hline
       & \multicolumn{3}{c|}{$\frac{h_{\rm f}}{h} = 1/8$}
       & \multicolumn{3}{c|}{$\frac{h_{\rm f}}{h} = 1/4$}
       & \multicolumn{3}{c|}{$\frac{h_{\rm f}}{h} = 1/2$} \\
$p_k$  & $a_1(p_k)$ & $b_1(p_k)$ & $b_2(p_k)$ 
       & $a_1(p_k)$ & $b_1(p_k)$ & $b_2(p_k)$ 
       & $a_1(p_k)$ & $b_1(p_k)$ & $b_2(p_k)$ \\ \hline
0.00 & 0.76738 & 0.60579 & 0.42216 & 0.77051 & 0.61120 & 0.42389 & 0.77363 & 0.61953 & 0.45295 \\ 
0.04 & 0.83462 & 0.61172 & 0.44778 & 0.84224 & 0.61607 & 0.45470 & 0.87242 & 0.63691 & 0.47535 \\ 
0.08 & 0.82739 & 0.60701 & 0.45371 & 0.83470 & 0.61024 & 0.46291 & 0.86400 & 0.62988 & 0.48633 \\ 
0.12 & 0.81649 & 0.60264 & 0.45711 & 0.82285 & 0.60580 & 0.46643 & 0.84984 & 0.62610 & 0.48880 \\ 
0.16 & 0.80142 & 0.59934 & 0.45584 & 0.80744 & 0.60510 & 0.45995 & 0.83017 & 0.62289 & 0.48759 \\ 
0.20 & 0.78410 & 0.59769 & 0.44867 & 0.78861 & 0.60230 & 0.45500 & 0.80852 & 0.62596 & 0.47106 \\ 
0.24 & 0.76246 & 0.59063 & 0.44922 & 0.76533 & 0.59494 & 0.45598 & 0.78215 & 0.62213 & 0.46478 \\ 
0.28 & 0.73555 & 0.57859 & 0.45631 & 0.73659 & 0.58273 & 0.46306 & 0.74857 & 0.61036 & 0.47016 \\ 
0.32 & 0.70230 & 0.56192 & 0.46838 & 0.70107 & 0.56562 & 0.47540 & 0.70553 & 0.59107 & 0.48468 \\ 
0.36 & 0.66179 & 0.54059 & 0.48470 & 0.65752 & 0.54327 & 0.49266 & 0.65062 & 0.56369 & 0.50746 \\ 
0.40 & 0.61221 & 0.51377 & 0.50533 & 0.60360 & 0.51457 & 0.51511 & 0.57676 & 0.52412 & 0.54163 \\ 
\hline
\end{tabular}\\[1ex]
\caption{Coefficient values describing the optimized finite difference operators for the fd5 fine scale operator.}
\label{tab:coeff_optfd5}
\end{table}

\subsection{Optimized regular grid finite elements}

\begin{figure}[thf]
\centering
\includegraphics[height=44mm]{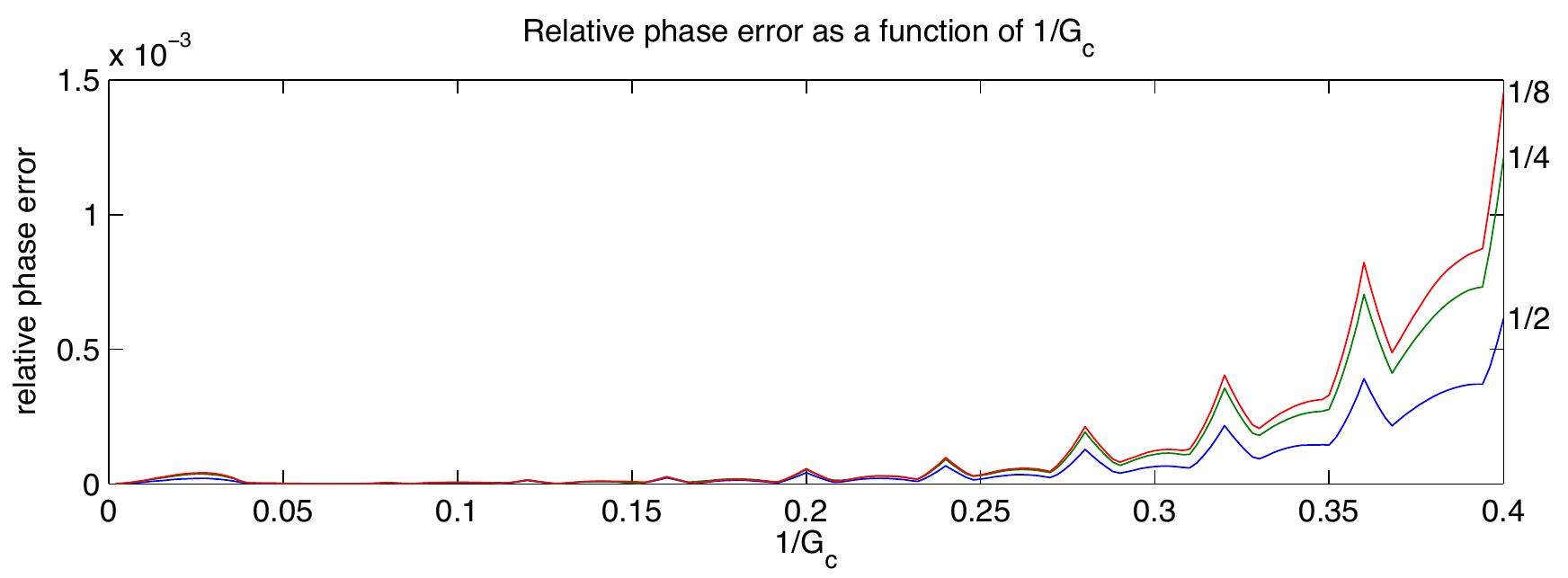}
\caption{Relative phase speed errors for the optimized scheme associated with (\ref{eq:discrete_op_FE}), maximum over $\theta$ as a function of $p = 1/G_{\rm c}$, for $h_{\rm f}/h = 1/2, 1/4$ and $1/8$.}
\label{fig:err_optFE}
\end{figure}

It is straightforward to derive the matrix for 2-D regular grid finite elements using first order rectangular (square) elements and constant $k^2$, see section~\ref{sec:further_numerics} below. One finds
\begin{equation} \label{eq:discrete_op_FE}
\begin{split}
  (H^{\rm FE}_h u)_{i,j} 
  = {}& 
( \frac{8}{3} - \frac{4}{9} h^2 k^2 ) u_{i,j}
+ ( - \frac{1}{3} - \frac{1}{9} h^2 k^2 ) ( u_{i+1,j} + u_{i-1,j} + u_{i,j+1} + u_{i,j-1} )
\\
{}& + ( - \frac{1}{3} - \frac{1}{36} h^2 k^2 )
 ( u_{i+1,j+1} + u_{i-1,j+1} + u_{i+1,j-1} + u_{i-1,j-1} ) .
\end{split}
\end{equation}
This equals (\ref{opt_fd_helm}) with $a_1 = 2/3, a_2=1/3, b_1=4/9, b_2=4/9$ and $b_3=1/9$, and an overal multiplicative factor $h^2$ due to the use of the finite element method instead of finite differences.
It is straightforward to find optimized coarse scale finite difference operators in the form
\begin{equation}
\begin{split}
    H^{\rm opt,FE}_h u (x_{i,j}) 
  = {}& (4 a_1 - h^2 k^2 b_1) u_{i,j} 
        + ( (-a_1+a_2) - h^2 k^2 b_2/4) (u_{i-1,j}+u_{i+1,j}+u_{i,j-1}+u_{i,j+1})\\
    {}& + ( -a_2 - h^2 k^2 b_3/4 ) (u_{i-1,j-1}+u_{i-1,j+1}+u_{i+1,j-1}+u_{i+1,j+1}) .
\end{split}
\end{equation}
The same optimization scheme is used as above. The values of the coefficients are given in Table~\ref{tab:coeff_optFE} and the relative phase speed differences as a function of $1/G_{\rm c}$ are given in Figure~\ref{fig:err_optFE}. 

\begin{table}
\begin{tabular}{|r|ccc|ccc|ccc|} \hline
       & \multicolumn{3}{c|}{$\frac{h_{\rm f}}{h} = 1/8$}
       & \multicolumn{3}{c|}{$\frac{h_{\rm f}}{h} = 1/4$}
       & \multicolumn{3}{c|}{$\frac{h_{\rm f}}{h} = 1/2$} \\
$p_k$  & $a_1(p_k)$ & $b_1(p_k)$ & $b_2(p_k)$ 
       & $a_1(p_k)$ & $b_1(p_k)$ & $b_2(p_k)$ 
       & $a_1(p_k)$ & $b_1(p_k)$ & $b_2(p_k)$ \\ \hline
0.00 & 0.76647 & 0.60220 & 0.42163 & 0.76485 & 0.59678 & 0.42071 & 0.75687 & 0.57482 & 0.41653 \\ 
0.04 & 0.82947 & 0.60801 & 0.44478 & 0.82174 & 0.60243 & 0.44031 & 0.79073 & 0.57900 & 0.42466 \\ 
0.08 & 0.82251 & 0.60347 & 0.45037 & 0.81512 & 0.59641 & 0.44891 & 0.78598 & 0.57656 & 0.42640 \\ 
0.12 & 0.81202 & 0.59934 & 0.45331 & 0.80514 & 0.59341 & 0.44959 & 0.77849 & 0.57352 & 0.42749 \\ 
0.16 & 0.79758 & 0.59608 & 0.45200 & 0.79216 & 0.59235 & 0.44396 & 0.76950 & 0.57382 & 0.41975 \\ 
0.20 & 0.78078 & 0.59376 & 0.44622 & 0.77726 & 0.59178 & 0.43482 & 0.75789 & 0.56950 & 0.41903 \\ 
0.24 & 0.75974 & 0.58635 & 0.44763 & 0.75885 & 0.58563 & 0.43396 & 0.74319 & 0.56157 & 0.42314 \\ 
0.28 & 0.73368 & 0.57436 & 0.45485 & 0.73543 & 0.57397 & 0.44094 & 0.72495 & 0.55058 & 0.43081 \\ 
0.32 & 0.70164 & 0.55808 & 0.46659 & 0.70647 & 0.55817 & 0.45237 & 0.70275 & 0.53691 & 0.44096 \\ 
0.36 & 0.66281 & 0.53747 & 0.48218 & 0.67130 & 0.53836 & 0.46744 & 0.67663 & 0.52099 & 0.45245 \\ 
0.40 & 0.61563 & 0.51184 & 0.50162 & 0.62872 & 0.51411 & 0.48582 & 0.64587 & 0.50233 & 0.46570 \\ 
\hline
\end{tabular}\\[1ex]
\caption{Coefficient values describing the optimized finite difference operators for the FE fine scale operator.}
\label{tab:coeff_optFE}
\end{table}

\subsection{Optimized finite difference scheme in 3-D}

In 3-D a class of operators like (\ref{opt_fd_helm}) can be constructed similarly. It depends on coefficients $a_1,a_2,a_3$ satisfying $a_1+a_2+a_3 = 1$ and on coefficeints $b_1,b_2,b_3,b_4$ satisfying
$b_1+b_2+b_3+b_4=1$, and is given by
\begin{equation} \label{opt_fd_helm3}
\begin{split}
& \hspace*{7mm}    H^{\rm opt}_h u (x_{i,j,k}) 
  = (6 a_1 h^{-2} - k^2 b_1) u_{i,j,k} 
\\
    {}& + ( (-a_1+a_2) h^{-2} - k^2 b_2/6) (u_{i-1,j,k}+u_{i+1,j,k}+u_{i,j-1,k}+u_{i,j+1,k}+u_{i,j,k-1}+u_{i,j,k+1})
\\
    {}& + ( (-2a_2 + a_3/2) h^{-2} - k^2 b_3/12 ) 
(  u_{i-1,j-1,k}+u_{i-1,j+1,k}+u_{i+1,j-1,k}+u_{i+1,j+1,k} + u_{i-1,j,k-1}
\\
{}& \;\;\;\;\;
             +u_{i-1,j,k+1}+u_{i+1,j,k-1}+u_{i+1,j,k+1}
 + u_{i,j-1,k-1}+u_{i,j-1,k+1}+u_{i,j+1,k-1}+u_{i,j+1,k+1} )
\\
    {}& + ( -3a_3/4 h^{-2} - k^2 b_4/8 ) 
(  u_{i-1,j-1,k-1}+u_{i-1,j+1,k-1}+u_{i+1,j-1,k-1}+u_{i+1,j+1,k-1}
\\
{}& \;\;\;\;\;
 + u_{i-1,j-1,k+1}+u_{i-1,j+1,k+1}+u_{i+1,j-1,k+1}+u_{i+1,j+1,k+1}) .
\end{split}
\end{equation}
There are now five coefficients that need to be determined by optimization. Again these are a function of $h_{\rm f}/h$ and of $p = 1/G_{\rm c}$. The procedure is done similarly as above. The vector $(\cos(\theta),\sin(\theta))$ is replaced by a unit vector on the sphere using spherical coordinates $(\theta,\phi)$, give by
$(\cos(\theta) \cos(\phi),\cos(\theta)\sin(\phi), \sin(\theta))$, and the objective functional is weighted by $\sin(\theta)$.

\begin{figure}[thf]
\centering
\includegraphics[height=44mm]{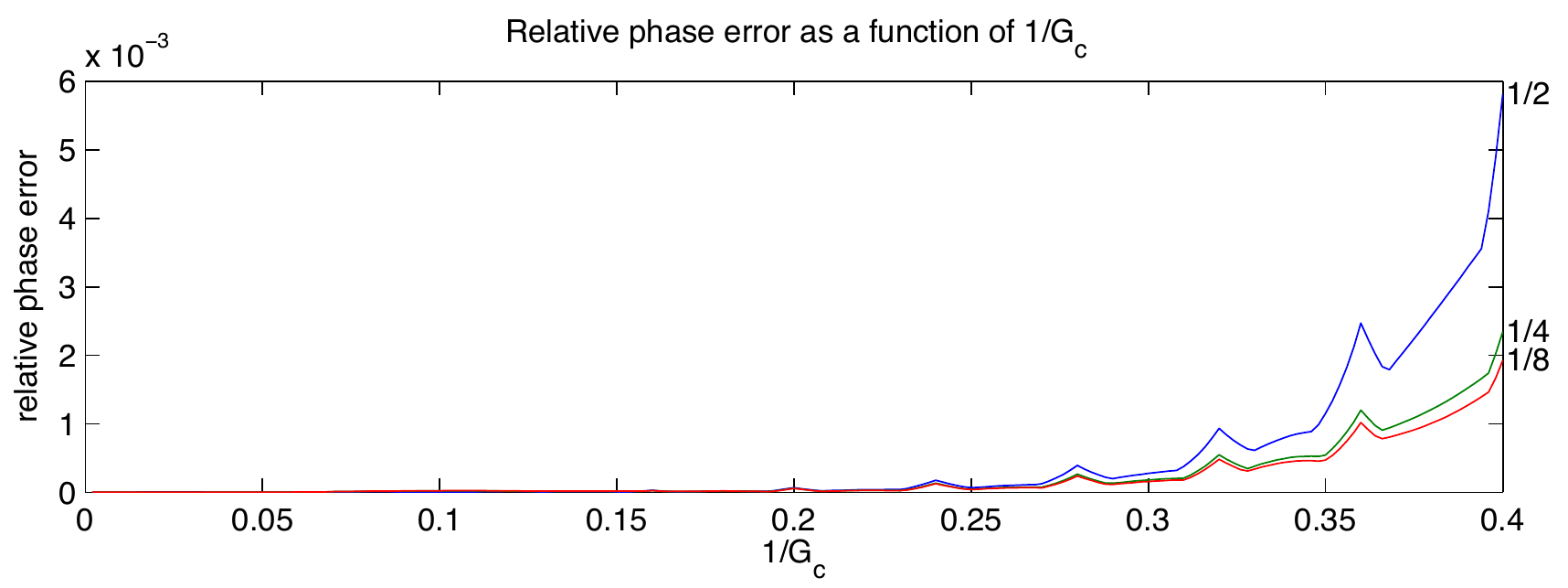}
\caption{Relative phase speed errors for the 3-D optimized scheme associated with the fd7 method, maximum over $(\theta,\phi)$ as a function of $p = 1/G_{\rm c}$, for $h_{\rm f}/h = 1/2, 1/4$ and $1/8$.}
\label{fig:err_optfd7}
\end{figure}

The same optimization scheme is used as above. The values of the coefficients are given in Table~\ref{tab:coeff_optfd7} and the relative phase speed differences as a function of $1/G_{\rm c}$ are given in Figure~\ref{fig:err_optfd7}. 

\begin{table}[t]
\begin{center}
\begin{tabular}{|r|ccccc|} \hline
       & \multicolumn{5}{c|}{$\frac{h_{\rm f}}{h} = 1/8$} \\
$p_k$  & $a_1(p_k)$ & $a_2(p_k)$ & $b_1(p_k)$ & $b_2(p_k)$ & $b_3(p_k)$ \\ \hline
0.00 & 0.75517 & 0.16259 & 0.54098 & 0.34422 & 0.19711 \\ 
0.04 & 0.75549 & 0.15831 & 0.54028 & 0.34418 & 0.19734 \\ 
0.08 & 0.74355 & 0.17283 & 0.54382 & 0.33897 & 0.19206 \\ 
0.12 & 0.70967 & 0.22754 & 0.54297 & 0.33616 & 0.19223 \\ 
0.16 & 0.69268 & 0.24313 & 0.54196 & 0.32594 & 0.20400 \\ 
0.20 & 0.67765 & 0.24848 & 0.53589 & 0.32075 & 0.21686 \\ 
0.24 & 0.65980 & 0.25228 & 0.52300 & 0.32683 & 0.22317 \\ 
0.28 & 0.63470 & 0.26238 & 0.50163 & 0.34886 & 0.21766 \\ 
0.32 & 0.60630 & 0.26831 & 0.48133 & 0.35550 & 0.23347 \\ 
0.36 & 0.58183 & 0.26501 & 0.46889 & 0.33974 & 0.26288 \\ 
0.40 & 0.55400 & 0.25602 & 0.45013 & 0.32489 & 0.29949 \\ 
\hline
       & \multicolumn{5}{c|}{$\frac{h_{\rm f}}{h} = 1/4$} \\ \hline
0.00 & 0.75957 & 0.16479 & 0.54568 & 0.34705 & 0.19853 \\ 
0.04 & 0.76194 & 0.16118 & 0.54572 & 0.34714 & 0.19860 \\ 
0.08 & 0.75567 & 0.16304 & 0.54554 & 0.34545 & 0.19740 \\ 
0.12 & 0.71958 & 0.22112 & 0.54645 & 0.33962 & 0.19835 \\ 
0.16 & 0.70279 & 0.23466 & 0.54490 & 0.33257 & 0.20536 \\ 
0.20 & 0.68712 & 0.23951 & 0.53872 & 0.33008 & 0.21307 \\ 
0.24 & 0.66708 & 0.24579 & 0.52468 & 0.34092 & 0.21302 \\ 
0.28 & 0.64032 & 0.25652 & 0.50408 & 0.36041 & 0.20978 \\ 
0.32 & 0.60968 & 0.26327 & 0.48374 & 0.36555 & 0.22774 \\ 
0.36 & 0.58186 & 0.26184 & 0.47007 & 0.35072 & 0.25739 \\ 
0.40 & 0.55039 & 0.25231 & 0.45013 & 0.33345 & 0.29951 \\ 
\hline
       & \multicolumn{5}{c|}{$\frac{h_{\rm f}}{h} = 1/2$} \\
0.00 & 0.77998 & 0.17505 & 0.56428 & 0.35970 & 0.20490 \\ 
0.04 & 0.78635 & 0.17442 & 0.56571 & 0.36071 & 0.20541 \\ 
0.08 & 0.78273 & 0.16881 & 0.56298 & 0.36150 & 0.20719 \\ 
0.12 & 0.76438 & 0.18678 & 0.56540 & 0.35620 & 0.20287 \\ 
0.16 & 0.74684 & 0.19603 & 0.56370 & 0.35299 & 0.20299 \\ 
0.20 & 0.72755 & 0.20131 & 0.55813 & 0.35277 & 0.20452 \\ 
0.24 & 0.70298 & 0.20847 & 0.54673 & 0.35830 & 0.20693 \\ 
0.28 & 0.66863 & 0.22424 & 0.52423 & 0.38368 & 0.19633 \\ 
0.32 & 0.62734 & 0.23845 & 0.49946 & 0.39740 & 0.20725 \\ 
0.36 & 0.58198 & 0.25329 & 0.47567 & 0.40216 & 0.22132 \\ 
0.40 & 0.53417 & 0.23589 & 0.45011 & 0.36784 & 0.29962 \\ 
\hline
\end{tabular}\\[1ex]
\end{center}
\caption{Coefficient values describing the 3-D optimized finite difference operators for the Fd7 fine scale operator.}
\label{tab:coeff_optfd7}
\end{table}

\section{The multigrid methods}
\label{sec:multigrid_methods}

The above described operator is used in two-grid and multigrid methods which we will now describe. For background on multigrid methods we refer to \cite{TrottenbergOosterleeSchueller2001}.  The multigrid method will be used as a {\em preconditioner} for a Krylov subspace method (GMRES \cite{SaadSchultz1986}).

We next discuss the choice of smoother. Smoothers apply one or more steps of an iterative method to the equation
\begin{equation}
  L_{h} u_{h} = f_h .
\end{equation}
This involves a choice of the iterative method and of the number of smoothing steps (presmoothing and postsmoothing). The result of a single iteration of the smoother will be denoted by $\operatorname{SMOOTH}(u_h,f_h)$. The $\omega$-Jacobi method and successive over relaxation (SOR) will be considered as smoothers. The latter is of course equal to Gauss Seidel (GS) when the relaxation factor is 1. 

The choices of grid refinement, and of the restriction and prolongation operators are standard. We start with a square grid with grid size $h$, and apply standard coarsening to square grids with grid size $2h$, $4h$ etc. 
For the restriction operator we always use full weighting, i.e.\ it is the tensor product of one-dimensional FW operators of the form,
in stencil notation,
$\frac{1}{4} \begin{bmatrix} 1 & 2 & 1 \end{bmatrix}$.
In 2-D it reads, in stencil notation
\begin{equation}
  R_h
  = \frac{1}{16} \begin{bmatrix}
    1 & 2 & 1 \\
    2 & 4 & 2 \\
    1 & 2 & 1
  \end{bmatrix} ,
\end{equation}
or, if $i,j$ are coordinates in the coarse grid
\begin{equation}
\begin{split}
  (R_h u)_{i,j} =  {}& 
  \frac{1}{16} \bigg(
  4 u_{2i,2j} + 2 ( u_{2i-1,2j} + u_{2i+1,2j} + u_{2i,2j-1} + u_{2i,2j+1} )
\\
  {}&
  + u_{2i-1,2j-1} + u_{2i+1,2j-1} + u_{2i-1,2j+1} + u_{2i+1,2j+1} \bigg)
\end{split}
\end{equation}
The prolongation operator is given by
\begin{equation}
  P_h = 4 ( R_{h} ) ' \, .
\end{equation}

For the two-grid method, the coarse grid correction operator is denoted by
\begin{equation}
  \operatorname{CGCORR}(u_h,f_h)
  = u_h + P_h ( L_{{\rm c},2h} )^{-1}  R_h ( f_h - L_h u_h ) .
\end{equation}
Here $L_h$ and $L_{{\rm c},2h}$ denote discretizations of $-\Delta -k^2$, that are still to be specified. 
We focus particularly on the choice of coarse grid operator $L_{{\rm c},2h}$. For the two-grid method in 2-D the following pairs of fine and coarse grid operators will be considered
\begin{enumerate}
\item
Both using standard 5-point finite differences (fd5-fd5)
\item
Using standard 5-point finite differences at the fine level and
a Galerkin approximation of this matrix at the coarse level (fd5-gal)
\item
Using standard 5-point finite differences at the fine level and the new optimized operator at the coarse level (fd5-opt)
\item
Using the operator from \cite{JoShinSuh1996} at the fine and coarse levels (jss-jss)
\end{enumerate}

The result of the two-grid cycle is obtained by first applying the smoother $\nu_1$ times, then the coarse grid correction and then the smoother $\nu_2$ times. Denoting by $u_{h,j}$, $j=1,\ldots,\nu_1+\nu_2$ the intermediate results, we define the map $\operatorname{TGCYC}(u_{h,0},f_h)$ by
\begin{equation}
  \operatorname{TGCYC}(u_{h,0},f_h) = u_{h,\nu_1+\nu_2+1} ,
\end{equation}
where
\begin{align}  
  u_{h,j+1}  = {}& \operatorname{SMOOTH}(u_{h,j},f_h) , 
        \qquad \text{ for $j=0,\ldots, \nu_1-1, \nu_1+1,\ldots \nu_1+\nu_2$}
\\
  u_{h,\nu_1+1} = {}& \operatorname{CGCORR}(u_{h,\nu_1},f_h) .
\end{align}
When used as a preconditioner it is applied with $u_{h,0} = 0$.

For multigrid with more than two coarsening levels we consider only the optimized method where at each scale the optimized operator is chosen to approximate the fine level finite difference operator. This is done using the parameter $h_{\rm f} / h$ that was introduced in section~\ref{sec:phase_speeds}. In this case the algorithm is based on the V-cycle. In the further developments most emphasis will be on the two-grid problem. Some experiments are done using multiple levels, to study whether the good two-grid performance extends to the multigrid method.

\section{Fourier analysis of multigrid methods}
\label{sec:MG_Fourier_analysis}

In this section, we compute and analyze two-grid convergence factors by using local Fourier analysis, cf.\ \cite{TrottenbergOosterleeSchueller2001}, chapters 3 and 4. We first describe their computation, and then give their values for a number of parameter choices. We identify choices of the smoother parameters that lead to small convergence factors. Implications of the divergence of the smoother for Helmholtz problems are briefly discussed.

As usual, the convergence of the following iteration is studied
\begin{equation} \label{eq:orig_MG_iteration}
  u^{(n+1)}_h = \operatorname{TGCYC}(u^{(n)}_h,f_h) .
\end{equation}
Let $u_{h,{\rm true}} = L_h^{-1} f_h$. The smoother acts linearly on the error $u_h - u_{h,\rm true}$, i.e.\ there is an operator $S$ such that
\begin{equation} \label{eq:smoothing_operator1}
  u_{h,{\rm true}} - \operatorname{SMOOTH}(u_h,f_h)
  = S_h ( u_{h,{\rm true}} - u_h ) .
\end{equation}
Similarly there is an operator $K_h^{2h}$ such that $\operatorname{CGCORR}(u_h,f_h)$ satisfies
\begin{equation}
  u_{h,{\rm true}} - \operatorname{CGCORR}(u_h,f_h)
  = K_h^{2h} ( u_{h,{\rm true}} - u_h ) . 
\end{equation}
For the $\operatorname{TGCYC}$ it follows that
\begin{equation} \label{eq:Moperator1}
  u_{h,{\rm true}} - \operatorname{TGCYC}(u_h,f_h)
  = M_h^{2h} ( u_{h,{\rm true}} - u_h ) . 
\end{equation}
with 
\begin{equation} \label{eq:operator_M}
  M_h^{2h} = S_h^{\nu_2} K_h^{2h} S_h^{\nu_1} . 
\end{equation}
In local Fourier analysis the spectral radius $\rho(M_h^{2h})$ is estimated. If $\rho(M_h^{2h}) < 1$ the method is asymptotically convergent.

Next we discuss the computation of the spectral radius. This is done in the dimensionless Fourier domain, i.e.\ the wave vector is written as $h^{-1} \theta = h^{-1} (\theta_1,\theta_2)$, so that the dimensionless wave numbers $\theta$ are in $[-\pi,\pi)^2$ (in 2-D). Sets of ``high'' and ``low'' wave numbers are defined by 
\begin{equation}
\begin{split}
  T^{\rm low} = {}& \left[- \frac{\pi}{2}, \frac{\pi}{2} \right)^2,
\\
  T^{\rm high} = {}& [-\pi,\pi)^2 \backslash T^{\rm low} . 
\end{split}
\end{equation}
For $\theta = (\theta_1, \theta_2) \in T^{\rm low}$ we also define
\begin{equation}
\begin{aligned}
  \theta^{(0, 0)} = {}& (\theta_1, \theta_2), \qquad\qquad
&
  \theta^{(1, 1)} = {}& (\bar \theta_1, \bar \theta_2) .
\\
  \theta^{(0, 1)} = {}& (\theta_1, \bar \theta_2),
&  
  \theta^{(1, 0)} = {}& (\bar \theta_1, \theta_2),
\end{aligned}
\end{equation}
where
\begin{equation}
\bar \theta_i = \left\{
                \begin{array}{cc}
                 \theta_i + \pi & \ \text{if } \ \theta_i <0\\
                 \theta_i - \pi & \ \text{if } \ \theta_i \ge 0,
                 \end{array}
                 \right.
\end{equation}
In local Fourier analysis, for each $\theta \in T^{\rm low}$ we define the four-dimensional space of harmonics by
\begin{equation}
  E^\theta = \operatorname{span} \big( \big\{ e^{i \theta^{\alpha} \cdot x / h} \, | \, 
    \alpha \in \{ (0,0),(1,1),(0,1),(1,0) \} \; \big\} \big)
\end{equation}
The operator $M_h^{2h}$ maps this space into itself, and hence can be reduced to $4 \times 4$ matrix for each $\theta \in T^{\rm low}$, which will be denoted by $\hat M_h^{2h}(\theta)$ and follows from (\ref{eq:operator_M}).

We next turn to the analysis of the smoothing operator. By $\hat{S}_h(\theta)$ we denote the $4 \times 4$ matrix that describes the action of $S_h$ on the space of harmonics. This matrix is diagonal, we write
\begin{equation}
  \hat{S}_h(\theta)
  = \operatorname{diag} (
  \tilde{S}_h(\theta^{(0,0)}),
  \tilde{S}_h(\theta^{(1,1)}),
  \tilde{S}_h(\theta^{(0,1)}),
  \tilde{S}_h(\theta^{(1,0)}) \, ) .
\end{equation}
In the $\omega$-Jacobi smoother the operator $L_{h}$ is written as the sum of its diagonal part $L_{{\rm diag},h}$ and offdiagonal part $L_{{\rm offdiag},h}$. The action of the $\omega$-Jacobi smoother is
\begin{equation}
  \operatorname{SMOOTH}(u_h,f_h)
  = (1-\omega_{\rm S}) u_h + \omega_{\rm S} ( L_{{\rm diag},h} )^{-1} ( f - L_{{\rm offdiag},h} u_h) .
\end{equation}
Using that $f_h = L_h u_{h,{\rm true}}$, it follows that $\tilde{S}_h(\theta)$ satisfies
\begin{equation}
  \tilde{S}_h(\theta) 
  = 1 - \omega_{\rm S} - \frac{\omega_{\rm S} \tilde{L}_{{\rm offdiag},h}(\theta)}{\tilde{L}_{{\rm diag},h}(\theta)}
\end{equation}
For the fd5 discretization we have $\tilde{L}_{{\rm diag},h} = 4 h^{-2} - k^2$, 
$\tilde{L}_{{\rm offdiag},h} = - 2 h^{-2} (\cos(\theta_1) + (\theta_2))$,
$\tilde{L}_{h} = h^{-2} (4 -2\cos(\theta_1) - 2\cos(\theta_2)) - k^2$.

The SOR smoother is given by
\begin{equation}
  \operatorname{SMOOTH}(u_h,f_h)
  = (1-\omega_{\rm S}) u_h + \omega_{\rm S} ( L_{+,h} )^{-1} ( f - L_{-,h} u_h) .
\end{equation}
where $L_{-,h}$ denotes the upperdiagonal part of $L_{-,h}$ and 
$L_{+,h} = L_{h} - L_{-,h}$. Similarly as above we obtain
\begin{equation}
  \tilde{S}_h(\theta) 
  = 1 - \omega_{\rm S} - \frac{\omega_{\rm S} \tilde{L}_{-,h}(\theta)}{\tilde{L}_{+,h}(\theta)}
\end{equation}
The explicit expressions for $L_{-,h}$ and $L_{+,h}$ are straightforward to derive. This concludes the Fourier domain analysis of the smoothers.

The Fourier transformed version of $K_h^{2h}$ is obtained in the standard way from the Fourier transformed versions of the operators that make up $K_h^{2h}$, i.e.
\begin{equation}
 \hat K_h^{2h}(\theta) = \hat I_h - \hat P_h (\theta) (\hat L_{{\rm c},2h}(2\theta))^{-1} 
    \hat R_h(\theta) \hat L_{h}(\theta) .
\end{equation}
Here $\hat I_h$ is the $4 \times 4$ identity matrix.
The operator $\hat L_{{\rm c},2h}(2 \theta)$ is a scalar given by the value of the symbol $\sigma(\theta^{(0,0)}/h)$. The operator $\hat L_h$ is a diagonal $4 \times 4$ matrix given by
\begin{equation}
  \hat{L}_h(\theta)
  = \operatorname{diag} (
  \tilde{L}_h(\theta^{(0,0)}),
  \tilde{L}_h(\theta^{(1,1)}),
  \tilde{L}_h(\theta^{(0,1)}),
  \tilde{L}_h(\theta^{(1,0)}) \, ) ,
\end{equation}
in which $\tilde{L}_h(\theta)$ is the scalar value of the symbol $\sigma(\theta / h)$. The Fourier transformed restriction and prolongation operators $\hat P_h$ and $\hat R_h$ are $4 \times 1$ and $1 \times 4$ matrices respectively, given by
\begin{equation}
 \hat{P}_h  = \begin{pmatrix}
        \frac{1}{4}(1+\cos      \theta_1)(1+\cos      \theta_2)\\
        \frac{1}{4}(1+\cos \bar \theta_1)(1+\cos \bar \theta_2)\\
        \frac{1}{4}(1+\cos \bar \theta_1)(1+\cos      \theta_2)\\
        \frac{1}{4}(1+\cos      \theta_1)(1+\cos \bar \theta_2) ,
      \end{pmatrix}
\end{equation}
and 
\begin{equation}
  \hat{R}_h = \hat{P}_h^T .
\end{equation}

The matrix $\hat{M}_h^{2h}(\theta)$ now follows from $\hat{S}_h(\theta)$ and $\hat{K}_h^{2h}(\theta)$ using (\ref{eq:operator_M}). The spectral radius $\rho_{\rm loc}(M_h^{2h})$ is given by
\begin{equation} \label{eq:rho_loc_M_all}
  \rho_{\rm loc} (M_h^{2h}) 
  = \sup \left\{ \rho_{\rm loc}(\hat M_h^{2h}(\theta)) \,|\, \theta \in T^{\rm low} \right\},
\end{equation}
and that $\rho_{\rm loc}(\hat M_h^{2h}(\theta))$ is the spectral radius of the $4 \times 4$ matrix $\hat M_h^{2h}(\theta)$.

\subsection{Divergence of the smoother}%
\label{subsec:smoother_divergence}

The fact that the smoother is divergent is easily shown using the smoothing factors. For example, with standard 5-point finite differences for the Helmholtz operator and $\omega$-Jacobi these are given by
\begin{equation} \label{eq:smoothing_factor_omega_jac}
  \tilde S_h(\theta) 
  = 1-\omega_{\rm S} + \frac{2\omega_{\rm S}}{4-k^2 h^2} (\cos \theta_1 + \cos \theta_2).
\end{equation}
For $\theta_1,\theta_2$ near zero this factor becomes $>1$, leading to the instability. For SOR something similar happens. This means that the coarse grid correction has to make up for the error introduced, and that performance in general is degraded. It also implies that one must be careful when increasing the number of smoothing steps.

\subsection{Computation of two-grid convergence factors}

The two-grid convergence factors are computed using the expression (\ref{eq:rho_loc_M_all}), taking the supremum for $\theta$ on a grid in polar coordinates. The grid was taken large enough to accurately compute the maximum in (\ref{eq:rho_loc_M_all}). We have computed the convergence factor for the four schemes listed in Section~\ref{sec:multigrid_methods}, i.e.\ the optimized schemes fd5-opt, jss-jss and the non-optimized schemes fd5-fd5 and fd5-gal. We first vary the smoother parameters to find a good smoother. Then we present the convergence factors as a function of damping and $G_{\rm c}$.

In Table~\ref{tab:conv_findS_opt} we present $\rho_{\rm loc} (M_h^{2h})$ for the fd5-opt method for two choices of smoother, $\omega$-Jacobi with $\omega_{\rm S} \in \{ 0.6,0.7,0.8,0.9, 1.0\}$ and SOR with $\omega_{\rm S} \in \{ 0.8, 0.9, 1.0, 1.1,1.2\}$, and six choices of the number of pre- and postsmoothing steps $(\nu_1,\nu_2)$, ranging from $(1,1)$ to $(6,6)$. The parameter $\alpha$ was fixed at $0.0025$, while $G_{\rm c}$ was chosen constant equal to 3.5. The numbers show that the performance of the scheme depends strongly on the parameters of the iterative method. Best performance for the optimized scheme occurs roughly for $\omega$-Jacobi with $\omega_{\rm S} =0.8$ and $(\nu_1, \nu_2) = (4, 4)$. 
For the jss-jss method a similar computation was made for $G_{\rm c} = 4$, and it was determined that $\omega$-Jacobi with $(\nu_1, \nu_2) = (2, 2)$ and $\omega_{\rm S}=0.8$ was either optimal or very close to optimal.
Table~\ref{tab:conv_findS_fd5} contains the results of a similar computation for the fd5-fd5 scheme. Now three choices of $(\nu_1,\nu_2)$ were used, $\alpha$ equaled $0.01$ or $0.02$ and $G_c$ was chosen equal to 10, because for lower values convergence was poor or absent. For the fd5-fd5 scheme the radius of convergence is less sensitive to the choice of smoother parameters  once $\alpha$ and $G_{\rm c}$ are such that reasonable  convergence can be obtained. We concentrate on the method with $(\nu_1, \nu_2) = (2, 2)$ and $\omega_{\rm S}=0.8$. Similarly, for fd5-gal the method is less sensitive to the detailed choice of parameters and we adopt the same choice.

In Table~\ref{tab:conv_alphaG} we compute the convergence radius as a function of $\alpha$ and $G_{\rm c}$, including results for the fd5-gal and jss-jss methods. We clearly observe that the fd5-opt method gives the best results, somewhat better than the jss-jss method and much better than fd5-fd5 and fd5-gal. It allows for the coarsest grids with the least amount of damping.

\begin {table}
\begin{center}
\small%
\begin{tabular}{|r||cccccc|}
\hline
$\alpha$& \multicolumn{6}{c|}{$0.0025$} \\ \hline
$(\nu_1, \nu_2)$  
      & $(1,1)$& $(2,2)$& $(3,3)$& $(4,4)$& $(5,5)$& $(6,6)$ \\ \hline
0.6J  &   $>1$ &   $>1$ &   $>1$ &  0.557 &  0.304 &  0.214  \\
0.7J  &   $>1$ &   $>1$ &  0.685 &  0.307 &  0.206 &  0.214  \\
0.8J  &   $>1$ &   $>1$ &  0.362 &  0.209 &  0.214 &  0.246  \\
0.9J  &   $>1$ &   $>1$ &   $>1$ &   $>1$ &   $>1$ &   $>1$  \\
   J  &   $>1$ &   $>1$ &   $>1$ &   $>1$ &   $>1$ &   $>1$  \\
\hline
SOR0.8&   $>1$ &   $>1$ &  0.426 &  0.393 &  0.729 &   $>1$ \\
SOR0.9&   $>1$ &  0.809 &  0.324 &  0.508 &   $>1$ &   $>1$ \\
    GS&   $>1$ &  0.527 &  0.321 &  0.657 &   $>1$ &   $>1$ \\
SOR1.1&   $>1$ &  0.476 &  0.381 &  0.846 &   $>1$ &   $>1$ \\
SOR1.2&   $>1$ &  0.452 &  0.453 &   $>1$ &   $>1$ &   $>1$ \\
\hline
\end{tabular}
\end{center}
\caption{Two-grid convergence factors $\rho_{\rm loc} (M_h^{2h})$, $G_{\rm c}=3.5$, fd5-opt.}
\label{tab:conv_findS_opt}
\end{table}

\begin{table}
\begin{center}
\small
\begin{tabular}{|r||ccc|ccc|}
\hline
  $\alpha$& & $0.01$  &&&$0.02$  &\\ \hline
$(\nu_1, \nu_2)$  
      & $(1,1)$& $(2,2)$& $(3,3)$& $(1,1)$& $(2,2)$& $(3,3)$\\  \hline
0.6J  &   $>1$ &   $>1$ &   $>1$ &  0.763 &  0.635 &  0.619 \\
0.7J  &   $>1$ &   $>1$ &   $>1$ &  0.697 &  0.622 &  0.617 \\
0.8J  &   $>1$ &   $>1$ &   $>1$ &  0.659 &  0.618 &  0.616 \\
0.9J  &   $>1$ &   $>1$ &   $>1$ &  0.677 &  0.617 &  0.616 \\
   J  &   $>1$ &   $>1$ &   $>1$ &   $>1$ &   $>1$ &   $>1$ \\
\hline
SOR0.8&   $>1$ &   $>1$ &   $>1$ &  0.621 &  0.615 &  0.615  \\
SOR0.9&   $>1$ &   $>1$ &   $>1$ &  0.604 &  0.616 &  0.615  \\
GS    &   $>1$ &   $>1$ &   $>1$ &  0.611 &  0.616 &  0.615  \\
SOR1.1&   $>1$ &   $>1$ &   $>1$ &  0.630 &  0.616 &  0.615  \\
SOR1.2&   $>1$ &   $>1$ &   $>1$ &  0.659 &  0.617 &  0.614  \\
\hline
\end{tabular}
\end{center}
\caption{Two-grid convergence factors $\rho_{\rm loc} (M_h^{2h})$, $G_{\rm c}=10$, fd5-fd5.}
\label{tab:conv_findS_fd5}
\end{table}

\begin {table}
\begin{center}
\small
\begin{tabular}{|r|ccc||ccc|} \hline
\multirow{2}{*}{$G_{\rm c}$}
       & \multicolumn{3}{c||}{fd5-opt}      &\multicolumn{3}{c|}{jss-jss}\\             
       & \hspace*{-1ex} $\alpha$=1.25e-3 \hspace*{-1ex} 
                          & 0.005 &  0.02 & \hspace*{-1ex} $\alpha$=1.25e-3 \hspace*{-1ex}
                                                             & 0.005 &  0.02 \\
\hline
3      &   0.634          & 0.439 & 0.438 &    $>1$          &  $>1$ &  $>1$ \\
3.5    &   0.228          & 0.204 & 0.202 &    $>1$          &  $>1$ & 0.502 \\
4      &   0.170          & 0.156 & 0.154 &    $>1$          & 0.639 & 0.231 \\
5      &   0.113          & 0.100 & 0.099 &    $>1$          & 0.293 & 0.122 \\
6      &   0.079          & 0.079 & 0.079 &    $>1$          & 0.358 & 0.121 \\
7      &   0.071          & 0.071 & 0.071 &    $>1$          & 0.359 & 0.110 \\
8      &   0.067          & 0.067 & 0.067 &    $>1$          & 0.324 & 0.095 \\
\hline\hline
\multirow{2}{*}{$G_{\rm c}$}
       & \multicolumn{3}{c||}{fd5-fd5}      &\multicolumn{3}{c|}{fd5-gal}\\             
       & \hspace*{-1ex} $\alpha$=1.25e-3 \hspace*{-1ex} 
                          & 0.005 &  0.02 & \hspace*{-1ex} $\alpha$=1.25e-3 \hspace*{-1ex}
                                                             & 0.005 &  0.02 \\
\hline
6      &    $>1$          &  $>1$ &  $>1$ &    $>1$          &  $>1$ &  $>1$ \\
7      &    $>1$          &  $>1$ &  $>1$ &    $>1$          &  $>1$ &  $>1$ \\
8      &    $>1$          &  $>1$ & 0.963 &    $>1$          &  $>1$ & 0.896 \\
10     &    $>1$          &  $>1$ & 0.618 &    $>1$          &  $>1$ & 0.588 \\
12     &    $>1$          &  $>1$ & 0.430 &    $>1$          &  $>1$ & 0.415 \\
\hline
\end{tabular}
\end{center}
\caption{Two-grid convergence factors as a function of $\alpha$, $G_{\rm c}$,
for fd5-opt, jss-jss, fd5-fd5 and fd5-gal.}
\label{tab:conv_alphaG}
\end{table}

\section{Numerical results}
\label{sec:numerical_results}

In this section we compare the convergence of multigrid methods with different choices of the coarse scale operators. We start with several experiments for the two-grid method in 2-D.
Then we investigate the multigrid method, the 3-D case and the case with PML boundary layers, in order to establish that the behavior observed in the 2-D, two-grid experiments also occurs more generally. 
For a full comparision we study the trade-off between damping present in the equation, the number of grid points per wave length at the coarse level $G_{\rm c}$ and the number of iterations required to reduce the residual by a factor of $10^{-6}$. 
In 2-D it is relatively easy to perform this kind of experiments. In 3-D the size of the coarse scale problem quickly becomes large for modern desktop machine (a machine with 8GB memory was used for these experiments using a Matlab implementation).

\subsection{The two-grid method in 2-D}
In the 2-D, two-grid examples we compare four methods. The regular five point operator at the fine scale is combined at the coarse scale with three choices of coarse scale operator, namely the regular five point operator, the Galerkin operator and our optimized method. These combinations are denoted fd5-fd5, fd5-gal and fd5-opt. We also study the JSS operator, used both at the fine and at the coarse scale, referred to as jss-jss. As shown below, sharp differences are present between the optimized coarse scale methods on the one hand and the fd5 and Galerkin coarse scale methods on the other hand.

The first series of experiments concerned a constant coefficient medium, i.e.\ $k$ is constant equal to $\frac{\pi}{G_{\rm c} h}$, where $G_{\rm c}$ is the number of grid points per wavelength at the coarse scale.
Experiments were performed with the $\omega$-Jacobi smoother with $\omega_{\rm S} =0.8$ and $\nu_1=\nu_2=4$ for the fd5-opt method, and $\omega_{\rm S}=0.8$ and $\nu_1 = \nu_2 = 2$ for the other methods.
The results are given in Table~\ref{tab:results_numexp1}. The following can be clearly observed
\begin{itemize}
\item
For the fd5-fd5 and fd5-gal method good convergence starts at $(\alpha,G_{\rm c}) = (0.01,10)$
or $(0.02,8)$. Relatively large $G_c$ and $\alpha$ are required.
\item
The jss-jss method perform well with
$G_{\rm c}$ as low as 3.5 and $\alpha = 0.0025$ or $\alpha = 0.005$ or larger.
\item 
In this example the fd5-opt method even performs well with $G_{\rm c} = 3$ and $\alpha = 0.00125$.
\end{itemize}
Here $N_{\rm iter} \lesssim 20$ is used as the (somewhat subjective) criterion for good convergence. Clearly the optimized methods (jss and opt) perform {\em much} better than the conventional choices fd5 or gal at the coarse level, with our fd5-opt method having the best performance.

Next we test this in two variable coefficient media. The first is a random medium, for which results are given in Table~\ref{tab:results_numexp2}, and the second is the Marmousi model, see Table~\ref{tab:results_numexp3}.
The Marmousi model has size 9200 x 3000 meters, and wave speeds between 1500 and 5500 ms${}^{-1}$, see the plot.
The value of $G_{\rm c}$ is the minimum value present in the model. The coefficients for a row of the matrix were obtained by freezing the value of $k$ locally. The performance of the opt-fd5 method for the Marmousi example is comparable to that in the constant medium case, the other methods perform somewhat better, probably because in most of the domain the actual value of $G_{\rm c}$ is larger than the minimum value. In these examples Dirichlet boundary conditions were used.

\begin{table}
\begin{center}
\begin{tabular}{|c|ccccc||ccccc|}
\hline&&&&&&&&&&\\[-1em]
\multirow{2}{*}{$G_{\rm c}$}
  & \multicolumn{5}{c||}{$N_{\rm iter}$ for fd5-fd5, 1023 $\times$ 1023} 
  & \multicolumn{5}{c|}{$N_{\rm iter}$ for fd5-gal, 1023 $\times$ 1023} \\ 
  & $\alpha$=1.25e-3 & 2.5e-3 & 0.005 & 0.01 & 0.02 
  & $\alpha$=1.25e-3 & 2.5e-3 & 0.005 & 0.01 & 0.02 \\ \hline&&&&&&&&&&\\[-1em]
     6 & $>$100 & $>$100 & $>$100 &     95 &     33 & $>$100 & $>$100 & $>$100 &     74 &     27\\
     7 & $>$100 & $>$100 & $>$100 &     55 &     22 & $>$100 & $>$100 & $>$100 &     47 &     20\\
     8 & $>$100 & $>$100 & $>$100 &     37 &     17 & $>$100 & $>$100 &     95 &     33 &     15\\
    10 & $>$100 & $>$100 &     51 &     21 &     11 & $>$100 & $>$100 &     47 &     20 &     11\\
    12 & $>$100 &     80 &     30 &     14 &      9 & $>$100 &     75 &     28 &     14 &      9\\
\hline \hline&&&&&&&&&&\\[-1em]
\multirow{2}{*}{$G_{\rm c}$}
  & \multicolumn{5}{c||}{$N_{\rm iter}$ for fd5-opt, 1023 $\times$ 1023} 
  & \multicolumn{5}{c|}{$N_{\rm iter}$ for jss-jss, 1023 $\times$ 1023} \\ 
  & $\alpha$=1.25e-3 & 2.5e-3 & 0.005 & 0.01 & 0.02 
  & $\alpha$=1.25e-3 & 2.5e-3 & 0.005 & 0.01 & 0.02 \\ \hline&&&&&&&&&&\\[-1em]
     3 &     15 &     11 &      8 &      7 &      7 & $>$100 & $>$100 & $>$100 &     51 &     22\\
   3.5 &      7 &      6 &      6 &      5 &      5 & $>$100 &     94 &     33 &     16 &     10\\
     4 &      6 &      5 &      5 &      5 &      5 &     49 &     21 &     12 &      8 &      7\\
     5 &      4 &      4 &      4 &      4 &      4 &     26 &     13 &      8 &      7 &      6\\
     6 &      4 &      4 &      4 &      4 &      4 &     37 &     15 &      9 &      7 &      6\\
     7 &      4 &      4 &      4 &      4 &      4 &     35 &     15 &      9 &      7 &      5\\
     8 &      4 &      4 &      4 &      4 &      4 &     26 &     13 &      8 &      6 &      5\\
\hline
\end{tabular}
\end{center}
\medskip
\caption{Convergence results for the two-grid method in a constant medium using different combinations of fine and coarse scale operators. The fine grid size is indicated.}
\label{tab:results_numexp1}
\end{table}

\newlength{\asize}
\newlength{\bsize}
\setlength{\asize}{76mm}
\setlength{\bsize}{76mm}
\begin{table}[p]
\begin{center}
\hspace*{-12mm}
\begin{minipage}[t]{\asize}
\begin{center}
(a) velocity\\
\includegraphics[width=0.9\bsize]{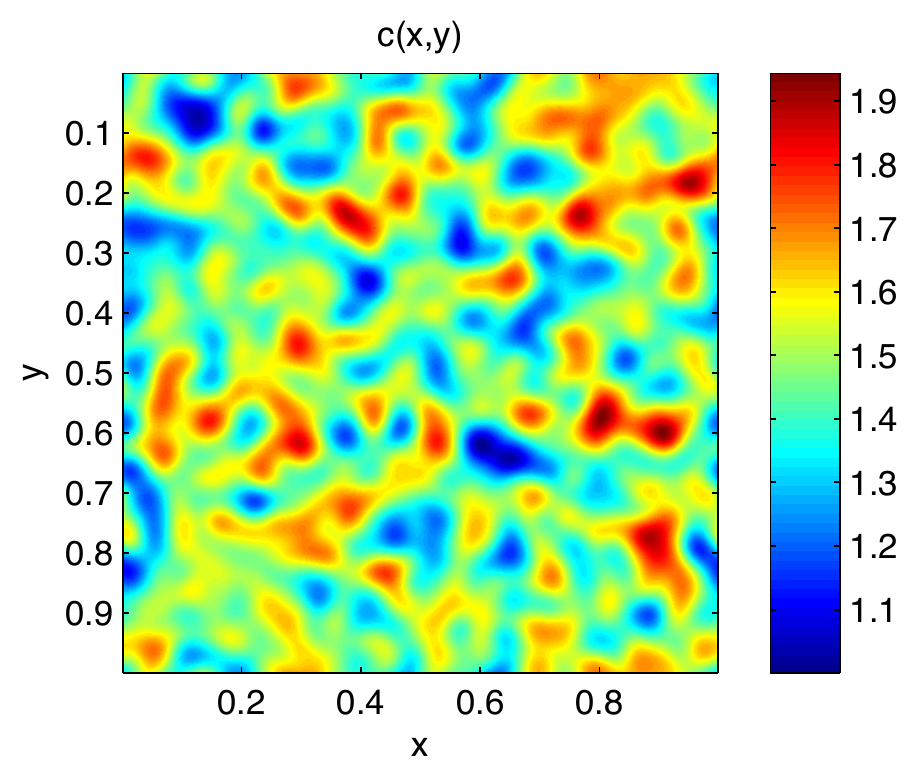}
\end{center}
\end{minipage}
\hspace*{2mm}
\begin{minipage}[t]{\bsize}
\begin{center}
(b) solution with $\text{nx}=511,G_{\rm c}=4,\alpha=\text{2.5e-3}$\\
\includegraphics[width=0.69\bsize]{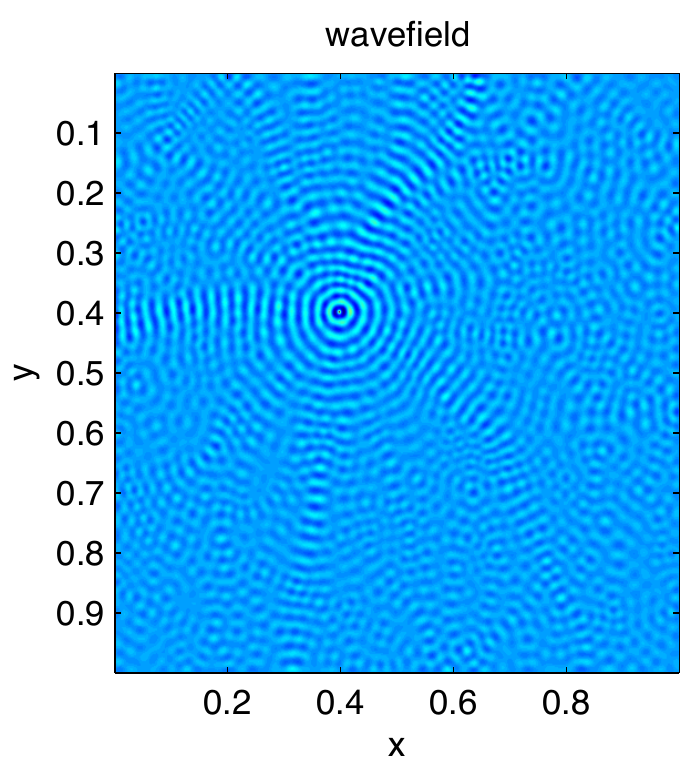}
\end{center}
\end{minipage}
\hspace*{-12mm}

\medskip
\begin{tabular}{|c|ccccc||ccccc|}
\hline&&&&&&&&&&\\[-1em]
\multirow{2}{*}{$G_{\rm c}$}
  & \multicolumn{5}{c||}{$N_{\rm iter}$ for fd5-fd5, 1023 $\times$ 1023} 
  & \multicolumn{5}{c|}{$N_{\rm iter}$ for fd5-gal, 1023 $\times$ 1023} \\ 
  & $\alpha$=1.25e-3 & 2.5e-3 & 0.005 & 0.01 & 0.02 
  & $\alpha$=1.25e-3 & 2.5e-3 & 0.005 & 0.01 & 0.02 \\ \hline&&&&&&&&&&\\[-1em]
     6 & $>$100 & $>$100 &     63 &     26 &     14 & $>$100 & $>$100 &     56 &     24 &     13\\
     7 & $>$100 & $>$100 &     39 &     18 &     11 & $>$100 & $>$100 &     36 &     17 &     10\\
     8 & $>$100 &     69 &     27 &     14 &      9 & $>$100 &     66 &     26 &     14 &      9\\
    10 & $>$100 &     36 &     17 &     10 &      7 &     97 &     34 &     16 &     10 &      7\\
    12 &     56 &     22 &     12 &      8 &      7 &     54 &     22 &     12 &      8 &      6\\
\hline \hline&&&&&&&&&&\\[-1em]
\multirow{2}{*}{$G_{\rm c}$}
  & \multicolumn{5}{c||}{$N_{\rm iter}$ for fd5-opt, 1023 $\times$ 1023} 
  & \multicolumn{5}{c|}{$N_{\rm iter}$ for jss-jss, 1023 $\times$ 1023} \\ 
  & $\alpha$=1.25e-3 & 2.5e-3 & 0.005 & 0.01 & 0.02 
  & $\alpha$=1.25e-3 & 2.5e-3 & 0.005 & 0.01 & 0.02 \\ \hline&&&&&&&&&&\\[-1em]
     3 &     5 &    5   &     5  &     5   &     5  &     25 &     15 &     11 &      8 &      7\\
   3.5 &     5 &    4   &     4  &     4   &     5  &     12 &      9 &      7 &      6 &      6\\
     4 &     4 &    4   &     4  &     4   &     4  &     12 &      9 &      7 &      6 &      6\\
     5 &     4 &    4   &     4  &     4   &     4  &     11 &      9 &      7 &      6 &      5\\
     6 &     4 &    4   &     4  &     4   &     4  &     10 &      8 &      6 &      5 &      5\\
     7 &     4 &    4   &     4  &     4   &     4  &      9 &      7 &      6 &      5 &      4\\
     8 &     4 &    4   &     4  &     4   &     4  &      8 &      6 &      5 &      5 &      4\\
\hline
\end{tabular}

\end{center}
\caption{Convergence results for different methods for a random medium. The fine grid size is indicated.}
\label{tab:results_numexp2}
\end{table}

%
%
\begin{table}[p]
\begin{center}
(a) medium\\
\includegraphics[width=104mm]{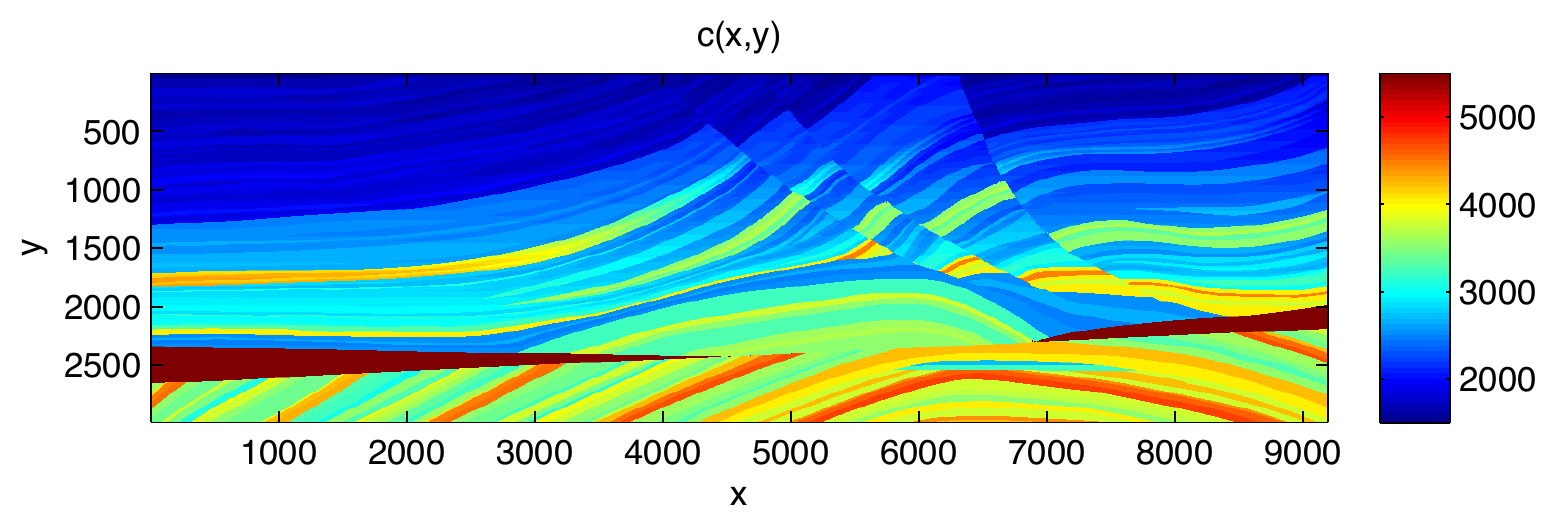}

\bigskip

(b) solution for $\text{nx}=1149, \alpha=\text{2.5e-3}, G_{\rm c}=4$\\
\includegraphics[width=94mm]{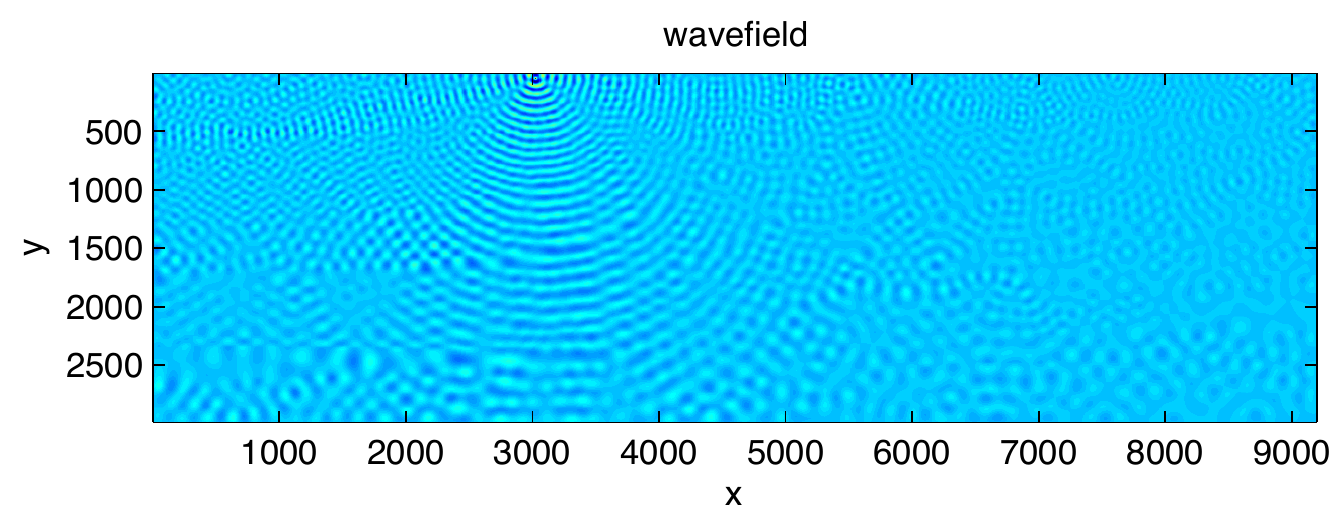}

\medskip

\begin{tabular}{|c|ccccc||ccccc|}
\hline&&&&&&&&&&\\[-1em]
\multirow{2}{*}{$G_{\rm c}$}
  & \multicolumn{5}{c||}{$N_{\rm iter}$ for fd5-fd5, 2299 $\times$ 749} 
  & \multicolumn{5}{c|}{$N_{\rm iter}$ for fd5-gal, 2299 $\times$ 749} \\ 
  & $\alpha$=1.25e-3 & 2.5e-3 & 0.005 & 0.01 & 0.02 
  & $\alpha$=1.25e-3 & 2.5e-3 & 0.005 & 0.01 & 0.02 \\ \hline&&&&&&&&&&\\[-1em]
     6 & $>$100 & $>$100 & $>$100 &     51 &     23 & $>$100 & $>$100 & $>$100 &     44 &     20\\
     7 & $>$100 & $>$100 &     70 &     30 &     16 & $>$100 & $>$100 &     63 &     27 &     15\\
     8 & $>$100 & $>$100 &     44 &     21 &     12 & $>$100 & $>$100 &     41 &     20 &     12\\
    10 & $>$100 &     52 &     24 &     14 &      9 & $>$100 &     50 &     23 &     13 &      9\\
    12 &     66 &     28 &     16 &     10 &      7 &     67 &     28 &     15 &     10 &      7\\
\hline \hline&&&&&&&&&&\\[-1em]
\multirow{2}{*}{$G_{\rm c}$}
  & \multicolumn{5}{c||}{$N_{\rm iter}$ for fd5-opt, 2299 $\times$ 749} 
  & \multicolumn{5}{c|}{$N_{\rm iter}$ for jss-jss, 2299 $\times$ 749} \\ 
  & $\alpha$=1.25e-3 & 2.5e-3 & 0.005 & 0.01 & 0.02 
  & $\alpha$=1.25e-3 & 2.5e-3 & 0.005 & 0.01 & 0.02 \\ \hline&&&&&&&&&&\\[-1em]
     3 &     17 &     13 &      9 &      7 &      6 & $>$100 & $>$100 &     75 &     29 &     16\\
   3.5 &     13 &     10 &      8 &      6 &      5 &     77 &     34 &     18 &     11 &      8\\
     4 &     11 &      9 &      7 &      5 &      5 &     19 &     12 &      9 &      7 &      7\\
     5 &      9 &      8 &      7 &      6 &      4 &     18 &     11 &      8 &      7 &      6\\
     6 &      8 &      7 &      6 &      5 &      4 &     16 &     10 &      8 &      6 &      5\\
     7 &      7 &      6 &      5 &      5 &      4 &     12 &      9 &      7 &      6 &      5\\
     8 &      7 &      6 &      5 &      5 &      4 &     11 &      8 &      7 &      5 &      5\\
\hline
\end{tabular}

\end{center}  

\caption{Convergence results for different methods for the Marmousi example. The fine grid size is indicated.}
\label{tab:results_numexp3}
\end{table}

\subsection{Further experiments: Multigrid, 3-D and PML layers\label{sec:further_numerics}}

Further numerical experiments are done to establish whether these results extend to multigrid experiments, to 3-D and to examples with PML layers. Our method to include PML layers is new and is described in this section.

Results for 2-D multigrid with 2, 3 and 4 levels, using a constant coefficient medium, are given in Table~\ref{tab:results_mg1}. In each case the same coarse level grid is used.
For most examples where the two-grid method converges reasonably fast, the multigrid method converges in about the same number of iterations. For the random medium in Table~\ref{tab:results_numexp2}, similar behavior is observed.

In 3-D the cost and memory use of the sparse factorizations scales worse than in 2-D. In this paper we study only a relatively small example of size $80^3$ that can be done on a regular machine with 8 GB memory using Matlab. This was approximately the largest example that could be done in this setup.
For these experiments we first had to determine suitable choices of the parameter $\omega_{\rm S}$ for $\omega$-Jacobi and for $\nu_1 = \nu_2$. This was done by test runs with $G_{\rm c} = 3.5$ and $\alpha=0.0025$. The best choices for $(\omega_{\rm S},\nu_1)$ were $\omega_{\rm S} = 0.8$ or $0.9$ and $\nu_1 = 7$ or $8$, we opted for $(\omega_{\rm S},\nu_1) = (0.9,8)$ which gave the fastest convergence. The results for the optimized method on the unit cube with constant wave speed are given in Table~\ref{tab:results_3d}. Generalizations of the method of \cite{JoShinSuh1996} to 3-D exist \cite{OpertoEtAl2007,ChenChengWu2012} but were not tested. In 3-D, the observed convergence behavior as a function of $\alpha$ and $G_{\rm c}$ does not differ much from the behavior observed in 2-D. Convergence is good for $G_{\rm c} \ge 3.5$.

We next consider the 2-D problem with PML boundary layers. 
This is important in practice -- rectangular domains are often encountered in combination
with PML boundary layers, for example in the seismic problem \cite{OpertoEtAl2007}.
As mentioned, we consider a finite element discretization for this problem.

As explained in \cite{Johnson_notespml}, in the PML method for the boundary $x_1 = \text{constant}$ the derivative $\frac{\partial}{\partial x_1}$ is replaced by 
$\frac{1}{1 + i \omega^{-1} \sigma_1(x_1)} \frac{\partial}{\partial x_1}$.
Here $\sigma_1$ is chosen to be 0 on the internal domain (no damping), and increases quadratically from the onset of the damping layer to boundary of the computational domain.
We introduce PML layers on all four boundaries.
It is convenient to multiply the equation by
$\left( 1 + i \omega^{-1} \sigma_1(x_1) \right)
\left( 1 + i \omega^{-1} \sigma_2(x_2) \right)$, this leads to a symmetric operator
easily discretized by finite elements.

\begin{figure}[p]
\begin{center}
\includegraphics[width=38mm]{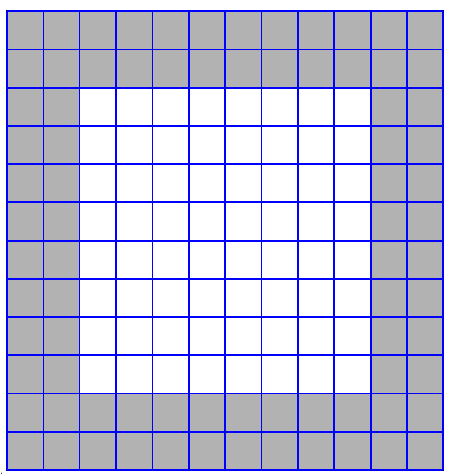}\hspace*{4mm}
\includegraphics[width=38mm]{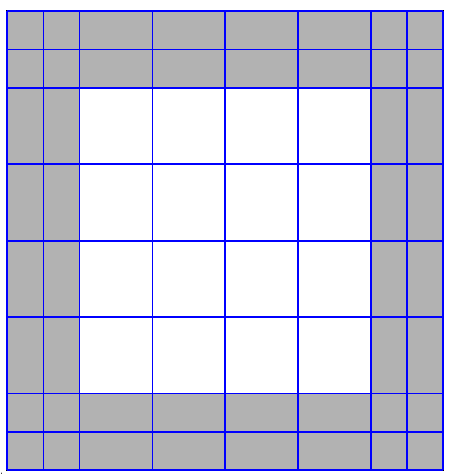}\hspace*{4mm}
\includegraphics[width=38mm]{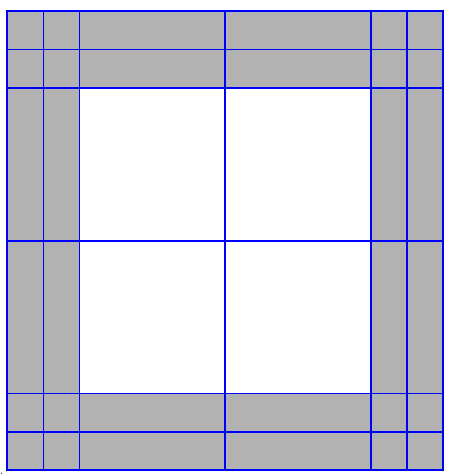}
\end{center}
\caption{Coarsening strategy to handle PML layers. Inside the PML layers (in grey) there is no coarsening in the direction normal to the boundary.}\label{fig:coarsening_PML}
\end{figure}

Straightforward inclusion of PML layers in the finite difference multigrid method
was observed to lead to very poor convergence, or no convergence at all.
Therefore we consider a coarsening strategy where inside the PML layer no coarsening takes place in the direction normal to the boundary. This leads to elongated basis functions in the boundary layers.
The resulting grids are given (schematically) in Figure~\ref{fig:coarsening_PML}.
We consider first order rectangular elements on these grids, denoted by $\phi_{ij}$.
The matrix elements are given by
\begin{equation} \label{eq:mat_elts_FE}
  M_{ij;kl} = \int \left[ \alpha_1 \alpha_2^{-1} \frac{\partial \phi_{ij}}{\partial x_1} 
        \frac{\partial \phi_{kl}}{\partial x_1} 
  + \alpha_1^{-1} \alpha_2 \frac{\partial \phi_{ij}}{\partial x_2} 
        \frac{\partial \phi_{ij}}{\partial x_2} 
  - k(x)^2 \alpha_1^{-1} \alpha_2^{-1} \phi_{kl}(x) \phi_{kl} \right] \, dx ,
\end{equation}
where the function $\alpha_j(x_j) = \frac{1}{1 + i \omega^{-1} \sigma_j(x_j)}$ contains the modifications for inclusion of the PML layer.

In the finite element context, the prolongation operator $P$ follows straightforwardly from specifying the grids, while the restriction operator $R$ is given by $P'$. 
Based on experiments with constant coefficients we choose as a smoother two iterations of the $\omega$-Jacobi method with relaxation factor $\omega_{\rm S} = 0.6$.

For the coarse grid operators the optimized stencils where computed in section~\ref{sec:phase_speeds}. These are used in the internal, square grid part of the domain, in such a way that the corresponding rows of the matrix contain the optimized stencil coefficients. For the other rows the regular finite element matrix elements are used.
We expect that the fact that no optimized coefficients are used inside the PML layer is of little importance, because the accuracy of the phase speeds is most relevant for long distance wave propagation, and there is no such propagation in the PML layer due to the damping. 

With this method a number of numerical experiments was carried out. Two choices of medium were considered, the constant and the random medium, both on the unit square with PML layers added outside it. A point source at $(0.3,0.3)$ was used a right hand side. We experimented with different multigrid levels from 2 to 4, using nx $\times$ ny= 200 $\times$ 200 grid points at the coarse level, and with different sizes for the two-grid method, using nx = ny = $200$, $400$ and $800$ at the coarse level. In each case values of $G_{\rm c}$ from 3 to 6 were used. The results are displayed in Table~\ref{tab:PML}. We see that the number of iterations is low for all cases with $G_{\rm c} \ge 3.5$.

Hence the new coarse grid operators can be used with PML boundary layers as well, providing similar reduction of $G_{\rm c}$ as for the examples without PML.

\begin{table}
\begin{center}
\begin{tabular}{|c|ccc||ccc||ccc|}
\hline&\multicolumn{9}{c|}{}\\[-1em]
\multirow{5}{*}{$G_{\rm c}$}
  & \multicolumn{9}{c|} {$N_{\rm iter}$ for fd5-opt, coarse size 255 $\times$ 255} \\ \hline&&&&&&&&&\\[-1em]
  & \multicolumn{3}{c||}{$\alpha=$1.25e-3} 
  & \multicolumn{3}{c||}{$\alpha=$0.005}
  & \multicolumn{3}{c|}{$\alpha=$0.02} \\ 
  & nlevels=2 & 3 & 4 
  & nlevels=2 & 3 & 4 
  & nlevels=2 & 3 & 4 \\ \hline&&&&&&&&&\\[-1em]
     3 &     15 &      8 &      8 &      8 &      7 &      8 &      7 &      7 &     8 \\
   3.5 &      7 &      6 &      6 &      6 &      5 &      6 &      5 &      6 &     6 \\
     4 &      6 &      5 &      5 &      5 &      5 &      5 &      5 &      5 &     5 \\
     5 &      4 &      4 &      4 &      4 &      4 &      4 &      4 &      4 &     4 \\
     6 &      4 &      4 &      4 &      4 &      4 &      4 &      4 &      4 &     4 \\
     7 &      4 &      4 &      4 &      4 &      4 &      4 &      4 &      4 &     4 \\
     8 &      4 &      4 &      4 &      4 &      4 &      4 &      4 &      4 &     4 \\
\hline
\hline
&\multicolumn{9}{c|}{}\\[-1em]
\multirow{5}{*}{$G_{\rm c}$}
  & \multicolumn{9}{c|} {$N_{\rm iter}$ for jss-jss, coarse size 255 $\times$ 255} \\ \hline&&&&&&&&&\\[-1em]
  & \multicolumn{3}{c||}{$\alpha=$1.25e-3} 
  & \multicolumn{3}{c||}{$\alpha=$0.005}
  & \multicolumn{3}{c|}{$\alpha=$0.02} \\ 
  & nlevels=2 & 3 & 4 
  & nlevels=2 & 3 & 4 
  & nlevels=2 & 3 & 4 \\ \hline&&&&&&&&&\\[-1em]
     3 & $>$100 & $>$100 & $>$100 & $>$100 & $>$100 & $>$100 &     22 &     22 &     22\\
   3.5 & $>$100 & $>$100 & $>$100 &     33 &     28 &     25 &     11 &     10 &     10\\
     4 &     46 &     25 &     28 &     11 &      9 &     10 &      7 &      8 &      8\\
     5 &     25 &     67 &     90 &      8 &     12 &     13 &      6 &      7 &      7\\
     6 &     33 &     61 &     74 &      9 &     12 &     12 &      6 &      6 &      6\\
     7 &     27 &     44 &     50 &      9 &     11 &     11 &      5 &      6 &      6\\
     8 &     18 &     24 &     26 &      8 &      9 &     10 &      5 &      6 &      6\\
\hline
\end{tabular}
\end{center}
\smallskip
\caption{Multigrid convergence for the constant medium with 2,3 and 4 levels, and
different values of $G_{\rm c}$ and $\alpha$}
\label{tab:results_mg1}
\end{table}

\begin{table}
\begin{center}
\begin{tabular}{|c|ccccc|}
\hline
\multirow{2}{*}{$G_{\rm c}$}
  & \multicolumn{5}{c||}{$N_{\rm iter}$ for fd7-opt} \\
  & $\alpha$=1.25e-3 & 2.5e-3 & 0.005 & 0.01 & 0.02 \\ \hline
     3 &     27 &     26 &     22 &     18 &     16 \\
   3.5 &      6 &      5 &      5 &      5 &      6 \\
     4 &      5 &      5 &      5 &      5 &      5 \\
     5 &      4 &      4 &      4 &      4 &      4 \\
     6 &      3 &      4 &      4 &      4 &      4 \\
     7 &      3 &      3 &      3 &      4 &      4 \\
     8 &      3 &      3 &      3 &      3 &      4 \\
\hline
\end{tabular}
\end{center}
\smallskip
\caption{Two-grid convergence for a constant coefficient medium in 3-D using the optimized coarse scale operator.}
\label{tab:results_3d}
\end{table}

\begin{table}
\begin{tabular}{|c|cccc|cccc|} \hline
  & \multicolumn{4}{c|}{constant medium} 
  & \multicolumn{4}{c|}{random medium} \\
  & $G_{\rm c}=3$ & 3.5 &   4 &   6 & $G_{\rm c}=3$ & 3.5 &   4 &   6 \\ \hline
nlevels=2, nxc=200 
  &  23 &  14 &  11 &   6 &  11 &   8 &   7 &   5 \\  
nlevels=3, nxc=200 
  &  24 &  14 &  11 &   6 &  10 &   8 &   7 &   5 \\
nlevels=4, nxc=200   
  & 100 &  17 &  12 &   8 &  25 &   8 &   7 &   6 \\ \hline\hline
nlevels=2, nxc=200
  &  23 &  14 &  11 &   6 &  11 &   9 &   7 &   5 \\
nlevels=2, nxc=400
  &  27 &  16 &  11 &   6 &  12 &   9 &   7 &   5 \\
nlevels=2, nxc=800
  &  35 &  17 &  12 &   6 &  13 &   9 &   7 &   5 \\ \hline
\end{tabular}
\smallskip
\caption{Iteration numbers for the numerical experiments with PML layer.}
\label{tab:PML}
\end{table}

\section{Conclusions and discussion}
\label{sec:discussion}

In this paper we have shown that the good or poor performance of the multigrid schemes for the Helmholtz equation is closely correlated with the phase speed differences between the fine and coarse scale operators. The results justify the conclusion that $G_{\rm c}$ can be reduced from about 10 to about 3.5, while at the same time reducing the amount of damping present, by using our new optimized finite differences as coarse scale operators.

The two-grid method now yields a general method to reduce the number of degrees of freedom in a high-frequency Helmholtz problem, at the cost of a few iterations. In high-frequency Helmholtz problems, often one makes use of optimized finite differences, or higher order methods at quite coarse grids, down to 10 or even less points per wavelength. With such coarse grids it is not obvious that multigrid can be applied. The new  two-grid method makes it possible to reduce the grid by a factor 2 in all directions once at least 7 points per wavelength are used.
One should also note that there is an obstruction against the direct use of the optimized discretizations at 3.5 points per wavelength or other schemes at this resolution, at least when reflections are present. The reason is that, according to linearized scattering theory, reflections are associated with Fourier components of the medium coefficient with wavevectors of length up to $2 \frac{\omega}{c}$, i.e.\ twice that of the wave field. The Nyquist criterion implies then that more than 4 points per wavelength are needed in the discretization of the medium coefficients in the fine scale operator. 
Our method is not restricted to the examples given. Optimized operators for other fine scale operators (or for the exact operator) can easily be constructed.

We have not compared the method with higher order (spectral) finite elements. Such a comparison should take into account the different requirements that methods may have concerning coefficient regularity, the sparsity patterns of the resulting matrices and the cost to invert them. This falls outside the scope of this paper.

We have deliberately opted for standard smoothers. This ensures that the observed good behavior is indeed due to the improved coarse scale operators. Other modifications, like polynomial smoothers for the PML method \cite{VanrooseRepsBinZubair2010,CalandraEtAl2013} were not needed here.
 
In 3-D the cost of the direct solver scales as $N^2$ for the sparse factorization and $N^{3/2}$ for the memory use \cite{George1973}. Our method can hence reduce the cost of the factorization by a factor $(10/3.5)^6 \approx 500$, and the memory use by a factor $(10/3.5)^{9/2} \approx 100$. However, this scaling also means that the coarse scale problem in general stays large. Further research should make clear whether the method can be used in combination with more efficient solvers, like the ones mentioned in the introduction.

\bibliographystyle{abbrv}
\bibliography{helmmgfd}

\end{document}